\documentclass[12pt]{article}

\usepackage{amscd, amsmath}
\usepackage{amsthm, amssymb}
\newcommand{\CC}{\mathbb{C}}
\newcommand{\RR}{\mathbb{R}}
\newcommand{\AAA}{\mathcal A}
\newcommand{\OO}{\mathcal O}
\newcommand{\si}{\sigma}
\newcommand{\sib}{\bar{\sigma}}
\newcommand{\db}{\bar{\partial}}
\newcommand{\dd}{{\rm d}}

\newcommand{\g}{\mathfrak g}
\newcommand{\ttt}{\mathfrak t}
\newcommand{\GG}{\mathcal G}
\newcommand{\UU}{\mathcal U}

\newcommand{\tvp}{\tilde{\varphi}}
\newcommand{\hh}{\mathfrak h}

\newtheorem{thm}{Theorem}[section]
\newtheorem{prop}[thm]{Proposition}
\newtheorem{cor}[thm]{Corollary}
\newtheorem{lem}[thm]{Lemma}

\newtheorem*{ques}{Question}

\theoremstyle{remark}
\newtheorem*{rem}{{\bf Remark}}

\theoremstyle{remark}
\newtheorem*{defn}{{\bf Definition}}

\theoremstyle{remark}
\newtheorem*{exa}{{\bf Example}}

\newenvironment*{prooff}{\noindent {\bf Proof.}}{\hfill $\qed$ \vspace{.3cm}}

\begin{document}

\begin{titlepage}
\title{
\vskip -70pt
\begin{flushright}
{\normalsize \ DAMTP-2004-64}\\
\end{flushright}
\vskip 35pt
{\bf Vortex equations in abelian \\ gauged $\sigma$-models}
}
\vspace{2cm}

\author{{J. M. Baptista} \thanks{ e-mail address:
    J.M.Baptista@damtp.cam.ac.uk }  \\
{\normalsize {\sl Department of Applied Mathematics and Theoretical
    Physics} \thanks{ address: Wilberforce Road, Cambridge
    CB3 0WA, England }
} \\
{\normalsize {\sl University of Cambridge}} 
}

\date{June 2004}

\maketitle

\thispagestyle{empty}
\vspace{1cm}
\vskip 20pt
{\centerline{{\large \bf{Abstract}}}}
\vspace{.35cm}
We consider nonlinear gauged $\si$-models with K\"ahler domain and
target. For a special choice of potential these models admit
Bogomolny (or self-duality) equations --- the so-called vortex
equations. Here we describe the space of solutions and energy spectrum
of the vortex equations when the gauge group is a torus $T^n$,
the domain is compact, and the target is $\CC^n$ or $\CC {\mathbb
  P}^n$. We also obtain a large family of solutions when the target is
a compact K\"ahler toric manifold. 

\end{titlepage}

\tableofcontents

\newpage

\section{Introduction}

Among the most general bosonic theories without gravity are the so-called
nonlinear gauged $\sigma$-models, also known as general Yang-Mills
theories with matter. These theories have been studied in the
theoretical physics literature for a long time now, and have recently
entered the mathematics literature as well.
To define them we need roughly the following data: two Riemannian
manifolds $M$ and $F$, a fibre bundle $E$ over the base $M$ with typical
fibre $F$, and a group $G$ acting on $F$ by isometries. The fields of the
theory are then a section $\phi : M \rightarrow E$ of the bundle and a
$G$-connection $A$. The energy functional is defined as
\begin{gather}
{\mathcal E}(A, \phi) \ =\ \int_M \ \frac{1}{2} \| F_A \|^2  +  \| \dd^A
\phi \|^2 + V(\phi) \ ,
\label{1.1}
\end{gather} 
where $F_A$ is the curvature of $A$, $\dd^A \phi $ is a covariant
derivative and $V(\phi)$ is a potential term. Notice that when the
bundle $E$ is trivial the section $\phi$ can be regarded as a map
$\phi : M\rightarrow F$, and so for $A=0$ this energy reduces to the
usual one for (non-gauged) $\sigma$-models.

In this paper we will be concerned with the case where $M$ and $F$ are
complex K\"ahler manifolds and the action of $G$ on $F$ is holomorphic
and hamiltonian. In this case there is a very special choice of
potential $V$, namely
\begin{gather}
V(\phi) \ =\ 2 \| \mu \circ \phi  \|^2 \ ,
\label{1.2}
\end{gather}
where $\mu$ is a moment map for the $G$-action on $F$. This
potential is special for two reasons, and it is a remarkable (though
not uncommon) fact that they occur simultaneously.

One reason is that with this choice the theory admits a supersymmetric
extension, at least when $M$ is an appropriate euclidean space. This
is an important fact well known in the physics literature (see for example
\cite{D-F}), but we will not make any use of it here. 
The other reason is that with the choice (\ref{1.2}) the energy
functional admits Bogomolny equations, or in other words has a
self-duality property. This fact appears to be less well known in
the physics literature, at least when the $\sigma$-model is nonlinear,
and apparently was first found in \cite{MiR, C-G-S}. When $M$ is a
Riemann surface these Bogomolny equations are 
\begin{align*}
& \db^A \phi \ =\ 0  \\
& \ast F_A + \mu \circ \phi \ =\  0 \ ,
\end{align*}
and can be generalized to any K\"ahler $M$. In this context these
equations are usually called vortex equations, because when $F=\CC$
and $G= U(1)$ they reduce to the usual vortex equations of the abelian
Higgs model. The solutions of these equations are exactly the global minima
(within each topological sector) of the energy functional ${\mathcal
 E}$. In fact they are also BPS states of the supersymmetric theory,
although we will not justify this here.
Hence it is usually interesting to know how many solutions these
equations admit up to gauge transformation, i.e. to describe the 
space of gauge equivalence classes of vortex solutions. For instance
in the abelian Higgs model ($F = \CC$ and $G= U(1)$) this was
originally done by Taubes for $M = \CC$ \cite{Tau} and by Bradlow for
$M$ compact K\"ahler \cite{Br}. In the more difficult non-abelian case
considerable progress has been made (e.g. \cite{Br1, Br-Da, Thad, Ba,
  MiR, H-T, E-I-N-O-S}), especially in the case where $F$ is a vector space, $G$
a unitary group and $M$ a Riemann surface. 

As a recent mathematical application, the vortex equations have been
used to define the so-called Hamiltonian Gromov-Witten
invariants \cite{MiR, C-G-S, C-G-MiR-S, MiR-T}. This is described from
a topological field theory point of view in \cite{Bap}.

$\ $

In this paper we will study the space of solutions of the
vortex equations for $M$ any compact K\"ahler manifold and $G$ an
abelian torus. At the end it turns out that we are able to completely describe 
this space in the case where $G=T^n$ and  $F=\CC^n$ or $F=\CC{\mathbb
  P}^n$. In some other cases a (big) family of non-trivial solutions 
is found, namely when $F$ is a compact K\"ahler toric manifold. 
The results obtained show an interesting interplay between the 
space of vortex solutions and the geometry of the moment polytope $\mu
(F)$ obtained from the torus action on $F$. An informal description of these
results is included in the comments of section 8. For the rest of this 
introduction we will just give a brief description of the content of each  
section.  
 
Section 2 is just a review of the model, where we try to carefully
describe all the notions involved in the definition of the energy
functional and of the equations. Since this nonlinear version of
gauge theory on fibre bundles with arbitrary fibres (as opposed to
vector space fibres) is not the most standard, we felt that this may
be useful. At the end of the section we also recall some standard
facts about complex gauge transformations and torus principal bundles,
which will be necessary further ahead.

In the short section 3 we give the space of solutions and energy spectrum of
the vortex equations in the case $G = T^n$ and $F = \CC^n$. When $n=1$
these are the classical vortex equations, defined on line bundles over
K\"ahler manifolds, and the solutions were described by Bradlow in
\cite{Br} and by Garc\'\i a-Prada in \cite{GP}. When $n>1$, following
work of Schroers \cite{Sch}, Yang has 
computed the space of solutions in the case where the
base is the complex plane or a compact Riemann surface \cite[p. 121]{Ya}.
The results contained in this section are for $n\ge 1$ and any compact K\"ahler
base. Their derivation follows quite straightforwardly from work in
\cite{Ba}. In the rest of the paper we will concentrate on the more
delicate case where $F$ is a compact manifold.

In section 4 we study the relation between the spaces of
solutions up to real gauge transformations and up to complex gauge
transformations. In fact, since the target $F$ is K\"ahler, the usual
$G$-gauge transformations can be extended to $G_\CC$-gauge
transformations. Then the first vortex equation is invariant under the
$G_\CC$-transformations whereas the second equation is invariant only
under the $G$-transformations. Thus it makes sense to ask if, given a
solution of the first equation, there exists a (unique) $G_\CC$-transformation
that takes it to a solution of the full equations. This question was
addressed by Mundet i Riera in \cite{MiR}, and a general ``stability''
criterion was 
found. This criterion, however, is generally not easy to evaluate in
practice. In section 4 we find that for $G=T^n$ and for suitable
conditions on $F$ this criterion is hugely simplified, and a direct
evaluation becomes possible. In particular, when $F$ is K\"ahler
toric, 
%%or whenever the polytope $\mu (F)  \subset \RR^n$ has $n+1$
%%vertices (and this includes the case $G=U(1)$), 
the answer to the question is essentially yes, and so there is an 
almost perfect correspondence between the real and 
complex moduli problems. The precise results are stated in section 4.1. 

In section 5 we determine the space of solutions and energy spectrum
of the vortex equations for $G=T^n$, $F=\CC {\mathbb P}^n$ 
and any compact K\"ahler $M$. The results obtained generalize the ones in
\cite{MiR} and 
\cite{S-S-Y}, where the authors determine same
quantities in the case where $M$ is a Riemann surface and $n=1$. The
calculations in this section require the results of section
4. The main results are stated in 5.1 and the proofs are contained in
5.2 and 5.3.

Section 6 is mainly preparatory. We study some general properties of
the vortex equations under quotients of the target
manifold $F$. Although we deal with a general group $G$, the results
will be mainly applied to $G=T^n$.

In section 7 we use the results of sections 4 and 6 to find
non-trivial solutions of the vortex equations for $G=T^n$ and $F$ a
compact K\"ahler toric manifold. This family is big enough so that
when $F=\CC {\mathbb P}^n$ it coincides with the full space of
solutions calculated in section 5. It is therefore natural to ask if,
for the other compact toric $F$, the solutions exhibited in this
section also exhaust the set of vortex solutions.

Finally in section 8 we make a few informal comments about the results
obtained. It may be helpful for the interested reader to have a look
at those before delving into the technicalities of the theorems.

\section{Review of the model}

\subsection{The energy functional}

The data we need to define the $\si$-model are the following.
\begin{itemize}
\item[$\bullet$] Two K\"ahler manifolds $M$ and $F$, with respective
K\"ahler forms $\omega_M$ and $\omega_F$.
\item[$\bullet$] A connected compact Lie group $G$, with Lie algebra $\g$,
and an Ad-invariant positive-definite inner product $\langle \ ,\
\rangle$ on $\g$.
\item[$\bullet$] An effective hamiltonian left action $\rho$ of $G$ on
  $F$ 
such that, for every $g\in G$, the transformations $\rho_g : F \rightarrow
F$ are holomorphic, and a moment map for this action $\mu : F \rightarrow 
\g^\ast$.
\item[$\bullet$] A principal $G$-bundle $\pi_P : P \rightarrow M$.
\end{itemize}
We remark that, in the fullest generality, the complex structure on
$F$ need not be assumed integrable, but we will assume that here.
Using the elements above one can define the associated bundle $E =
P\times_\rho F$, which is a bundle over $M$ with typical fibre $F$. It is
defined as the quotient of $P\times F$ by the equivalence relation $(p,q)
\sim (p\cdot g , g^{-1}\cdot q)$, for all  $g\in G$. The bundle    
projection $\pi_E : E \rightarrow M$ is determined by $\pi_E \circ \chi (p,q)
= \pi_P (p)$, where $\chi : P\times F \rightarrow E$ is the quotient map.
As a matter of notation, we will sometimes denote the equivalence
class $\chi (p,q)$ simply by $[p,q]$.

\vspace{.2cm}

\begin{defn}
The convention used here is that a moment map for the action $\rho $
of $G$ on $(F, \omega_F )$ is a map $\mu : F \rightarrow \g^\ast $ such
that
\begin{itemize}
\item[(i)]$\dd\, (\mu , \xi ) = \iota_{\xi^\flat}\, \omega_F $ in $\Omega^1
  (F)$ for all $\xi \in \g $, where $\xi^\flat$ is the vector field on $F$
  defined by the flow $t \mapsto \rho_{\exp (t\xi )}$.
\item[(ii)] $\rho_g^\ast \mu = {\rm Ad}_g^\ast \circ \mu $ for all $g
  \in G$, where ${\rm Ad}_g^\ast $ is the coadjoint representation of
  $G$ on $\g^\ast$.
\end{itemize}
If a moment map $\mu$ exists, it is not in general unique, but all the
other moment maps are of the form $\mu + a$, where $a \in [\g ,\g ]^0
\subset \g^\ast $ is a constant in the annihilator of $[\g, \g]
$. Recall also that under the identification $\g^\ast \simeq \g$
provided by an ${\rm Ad}$-invariant inner product on $\g$, the annihilator
$[\g, \g ]^0 $ is taken to the centre of $\g$. 
\end{defn}

\vspace{.2cm}

The fields of the theory are a connection $A$ on the principal bundle $P$
and a smooth section $\phi $ of $E$. Calling $\AAA$ the space of such   
connections and $\Gamma (E)$ the space of such sections, we define the
energy functional $\mathcal{E} : \AAA \times \Gamma (E) \rightarrow    
\RR_{\ge 0}$ of the $\si$-model by
\begin{gather}
\mathcal{E}(A,\phi )\ =\ \int_M \;  \left\{  \frac{1}{a^2} \| F_A \|^2  +
\| \dd^A   \phi \|^2 + a^2 \| \mu \circ \phi \|^2  \right\} \
\omega_M^{[m]}  \; ,\qquad  a\in \RR_{>0} .
\label{2.1}
\end{gather}
In this formula, as throughout the paper, $m$ is the complex dimension
of $M$, and we use the notation $ \omega_M^{[k]}  :=  \omega_M^k  /
k!$  for any $ k \in {\mathbb N}$. In particular $\omega_M^{[m]}$ is the metric
volume form on $M$. 

The various terms under the integral sign have the following meaning. 
$F_A$ is the curvature of the connection $A$. It can be regarded as a locally
defined $2$-form on $M$ with values in the Lie algebra $\g$. The norm   
$\| F_A \|^2 $ is then the natural one, induced simultaneously by the   
K\"ahler metric on $M$ and by  the inner product $\langle \ ,\ \rangle$ on
$\g$. 
In the third term of (\ref{2.1}), the norm $\| \cdot \| $ on $\g^\ast$ comes  
from the inner product  $\langle \ ,\ \rangle$, which induces an inner    
product on $\g^\ast$. By the ${\rm Ad}$-invariance of $\langle \ ,\ \rangle$  
and by the $G$-equivariance of the moment map $\mu$, the function $E
\rightarrow \RR_{\ge 0}$ determined by $\chi (p,q) \mapsto \|  \mu (q)\|^2$ is
well defined; the third term is then the composition of this function with
$\phi$.

As for the second term, its description is a little longer, since one
should first explain the meaning of the covariant derivative $\dd^A \phi$.
This is an extension of the usual notion of covariant derivatives on
vector bundles. We start by considering the differential of the quotient
map, $\dd \chi : TP \times TF \rightarrow TE$. A connection $A$ on  $P$
induces a horizontal distribution $H_A$ on $P$. Defining
$\mathcal{H}_A = \dd \chi (H_A)$, it is not difficult to  
show that the restrictions
\begin{gather}
\dd \pi_E : \mathcal{H}_A \longrightarrow TM  \qquad {\rm and} \qquad \dd   
\chi_{(p,q)} : T_q F \longrightarrow \ker(\dd \pi_E)_{\chi (p,q)}
\label{2.2}
\end{gather}
are isomorphisms, and in particular we get the splitting
\begin{gather}
TE \ =\ {\mathcal H}_A \oplus \ker \dd \pi_E \ .
\label{2.3}
\end{gather}
The covariant derivative of a section $\phi : M \rightarrow E$ is then
defined as the composition
\[
\begin{CD}
\dd^A \phi : TM @>{\dd \phi}>> TE = {\mathcal H}_A \oplus \ker \dd
\pi_E  @>{{\rm proj}_2}>>  \ker \dd \pi_E \ ,
\end{CD}
\]
where ${\rm proj}_2$ is just the projection. Notice that the image of $\dd^A  
\phi$ is in the tangent space to the fibres of $E$, which are isomorphic  
to $F$. Thus when $F$ is a vector space, the canonical isomorphism $T_v F 
\simeq F$ allows us to regard $\dd^A \phi$ as a map of vector bundles $TM 
\rightarrow E$, that is a section of $T^\ast M \otimes E$, which is the   
usual notion of covariant derivative on a vector bundle.
The norm $\| \dd^A \phi  \|^2 $ is defined in the usual way, using the
metric $g_M$ on  $M$ and the metric $g_F$ --- transported by the second
isomorphism of (\ref{2.2}) --- on $\ker \dd \pi_E $.

Finally notice that the constant $a^2$ can be absorbed by rescaling
the inner product on $\g$.

\subsection{The vortex equations}

Having explained the meaning of the energy functional (\ref{2.1}), we will now
see how to manipulate it in order to get Bogomolny equations.
First of all, using the isomorphisms (\ref{2.2}) and the splitting
(\ref{2.3}), one can 
transport the complex structures $J_M$ and $J_F$ of $M$ and $F$,
respectively, as well as the K\"ahler metrics $g_M$ and $g_F$, to the 
tangent bundle $TE$, thus defining a complex structure and a metric on
$TE$ by
\begin{gather}
J(A) = J_M \oplus J_F \qquad {\rm and} \qquad g(A) = g_M \oplus g_F \ .
\label{2.7}
\end{gather}
These depend on the connection $A$. Because the metrics $g_M$ and $g_F$
are K\"ahler, $J(A)$ is always compatible with $g(A)$, and so $(E , J(A),
g(A))$ is an almost-Hermitian manifold. Using this complex structure on $E$
and the one on $M$, one obtains a splitting $\dd^A \phi = \partial^A \phi
+ \db^A \phi$ by the usual formulae
\begin{align}
\db^A \phi\ &=\ \frac{1}{2} (\dd^A \phi + J_F\, \circ\, \dd^A \phi\,
\circ\, J_M )\ =\ 
\frac{1}{2}\, {\rm proj}_2 \,\circ \,(\dd \phi + J(A)\, \circ\, \dd \phi\, \circ
\, J_M )   \label{2.10} \\  
\partial^A \phi\ &=\ \frac{1}{2} (\dd^A \phi  - J_F\, \circ\, \dd^A\,
\phi\, \circ\, J_M ) \ .
\label{2.4}
\end{align}
For later convenience we also record here the local (i.e
trivialization-dependent) formulae for $\dd^A \phi$ and $\db^A
\phi$. Let $s : \UU \rightarrow P$ be a local section of $P$ over a
domain $\UU$ in $M$. Since $E=P\times_\rho F$ is an associated bundle,
this determines a trivialization of $E|_\UU$ by
\begin{gather}
\UU \times F \ \rightarrow \ E|_\UU \ , \qquad (x,q) \ \rightarrow \
[ s(x), q ] \ .
\label{2.11}
\end{gather}
With respect to these trivializations a section $\phi$ of $E$ can be
locally identified with a map $\hat{\phi} : \UU \rightarrow F$, and a
connection $A$ on $P$ can be identified with the connection form
$s^\ast A = \alpha\   \in \Omega^1 (\UU ; \g)$. Then the covariant
derivatives $\dd^A \phi$ and $\db^A \phi$ in $\Gamma (T^\ast M
\otimes\, \phi^\ast \ker \dd\, \pi_E )$ are locally given by 
\begin{align}
(\dd^A \phi)_q \ &=\ (\dd \, \hat{\phi})_q \ +\ (\alpha^l)_q \ \;
\xi_l^\flat |_{\hat{\phi } (q)}  \nonumber \\
(\db^A \phi)_q \ &=\ (\db \, \hat{\phi})_q \ +\ (\alpha^l)_q^{0,1} \ \;
\xi_l^\flat |_{\hat{\phi } (q)}  \qquad \qquad \forall\ q\in \UU \ ,
\label{2.12}
\end{align}
which are 1-forms on $T_q M$ with values in $T_{\hat{\phi} (q)} F$. In these
formulae $\{\xi_l \}$ is any basis for $\g$, the $\xi_l^\flat$ are the
vector fields on $F$ described in the definition of moment map
(section 2.1), and we have decomposed $\alpha = \alpha^l\, \xi_l$.

$\ $

We now come to the basic fact of the theory. This was first obtained in
\cite{MiR} and, for $M$ a Riemann surface, in \cite{C-G-S}.

\begin{thm}[\cite{MiR} , \cite{C-G-S}] 
For any connection $A \in  \AAA $ and any section $\phi \in \Gamma (E)$,
\begin{gather}
{\mathcal E} (A,\phi )\  =\ T_{[\phi ]} + \int_M \; \left\{ \| \frac{1}{a}
\Lambda F_A + a\ 
\mu\circ\phi  \|^2   +  2\| \db^A  \phi \|^2 +  \frac{4}{a^2} \| F_A^{0,2}
\|^2 \right\}\; \omega_M^{[m]} \ ,               
\label{2.5}
\end{gather}
where the term
\begin{gather}
T_{[\phi ]} \ =\  \int_M  \phi^{\ast} [\eta_E ] \wedge \omega_M^{[m-1]} -
\frac{1}{a^2} B_2 (F_A, F_A) \wedge \omega_M^{[m-2]}
\label{2.6}
\end{gather}
does not depend on $A$, and only on the homotopy class of $\phi $.
\label{t2.1}
\end{thm}

\begin{rem}
As is usually the case with these Bogomolny-type manipulations, there is
an alternative formula for ${\mathcal E}(A,\phi )$ which gives rise to the
anti-Bogomolny equations. This formula can be obtained from the one above 
by changing the sign of the first term of $T_{[\phi ]}$, substituting
$\db^A \phi$ for $\partial^A \phi$, and changing the plus to a minus sign
inside the first squared norm. The proof of \cite{MiR} is still
applicable, with minimal changes.
\end{rem}

\begin{cor}[\cite{MiR} ,\cite{C-G-S}] 
Within each homotopy class of the sections $\phi$ we have that
${\mathcal E}(A,\phi)
\geq  T_{[\phi ]}$, and there is an equality if and only if the pair
$(A,\phi )$ in $\AAA \times \Gamma (E)$ satisfies the equations
\begin{subequations}  \label{2.0}
\begin{align}
 &  \db^A \phi \ =\  0     \label{2.0a}   \\
 &  \Lambda F_A + a^2 \ \mu\circ\phi \ = \ 0   \label{2.0b}  \\
 &  F_A^{0,2}\ =\ 0  \ \ .  \label{2.0c} 
\end{align}
\end{subequations}
These first order equations are usually called vortex equations.
\end{cor}
Apart from $\db^A \phi$, several new terms appear in (\ref{2.5}) when compared
with (\ref{2.1}); their meaning is the following. The operator $\Lambda :   
\Omega^\ast (M) \rightarrow \Omega^{\ast -2}(M)$ is the adjoint, with   
respect to $g_M$, of the operator $\eta \mapsto \omega_M \wedge \eta $
on $\Omega^{\ast }(M)$. By well known formulae,
\begin{gather}
\Lambda F_A \ =\ \ast (\omega_M \wedge \ast F_A )\ =\  g_M (F_A ,
\omega_M)\ ,
\label{2.111}
\end{gather}
and so $\Lambda F_A$ can be seen as a locally defined function on $M$ with
values in $\g$, just as $\mu \circ \phi$. (More properly, they should be both regarded as global sections of $P\times_{{\rm   
Ad}_G}\g $.) Next, $F_A^{0,2}$ is just the $(0,2)$-component of $F_A$ under
the usual decomposition $\Omega^2 (M) = \Omega^{2,0}\oplus \Omega^{1,1}  
\oplus \Omega^{0,2}$. The form $B_2 (F_A , F_A)$ can be explicitly
written as 
\begin{gather}
B_2 (F_A , F_A) \ =\ F_A^j \wedge F_A^k \ \langle \xi_j , \xi_k \rangle \ ,
\label{2.8}
\end{gather}
where $\{ \xi_j \}$ is a basis of $\g$ and we have decomposed $F_A =
F_A^j \; \xi_j$; it represents the characteristic class of $P$
associated with the Ad-invariant polynomial $\langle 
\cdot , \cdot  \rangle : \g \times \g  \rightarrow \RR$.
 
Finally $[\eta_E]$ is a cohomology class in $H^2 (E)$, and is defined as  
follows. Consider the $2$-form on $P\times F$
\begin{gather}
\eta (A) \ := \ \omega_F - \dd (\mu , A) \ ,
\label{2.9}
\end{gather}
where we regard the connection $A$ as a form in $\Omega^1 (P,\g )$, in
the usual sense, and  $(\cdot , \cdot ) : \g^\ast \times \g \rightarrow
\RR $ is the natural pairing. The quotient map $\chi : P\times F
\rightarrow E$ is in a natural way a principal $G$-bundle, and it is not
difficult to check that the form $\eta (A)$ is invariant under the
associated $G$-action $(p,q) \mapsto (p\cdot g , g^{-1}\cdot q)$ on  
$P\times F$. Furthermore $\eta (A)$ is also a horizontal form, in the   
sense that it annihilates vectors in $\ker \dd \chi$, and so $\eta 
(A)$ descends to $E$, that is $\eta (A) = \chi^\ast \eta_E (A)$ for some
$\eta_E (A)$ in $\Omega^2 (E)$. The form $\eta_E (A)$ on $E$ is sometimes
called the minimal coupling form. Now, since $\eta (A)$ is closed,
$\eta_E (A)$ is also closed, and it is not difficult to show  that its
cohomology class in $H^2 (E)$ does not depend on $A$.
 We can therefore define $[\eta_E ]$ to be the cohomology class of
the forms $\eta_E (A)$.

\vspace{.2cm}

\begin{rem}
There is another way to look at the class $[\eta_E ]$ on $H^2 (E)$, using
the Cartan complex for the $G$-equivariant cohomology of  $F$. In this   
context, $[\eta_E ]$ is just the image by the Chern-Weil homomorphism of
the cohomology class in $H^2_G (F)$ determined by the equivariantly closed
form $\omega_F - X^b \mu_b  \ \in \ \Omega_G^2 (F)$ (see for example
\cite[ch. VII]{B-G-V}) .
\end{rem}

\vspace{.2cm}

Observe that the term $T_{[\phi]}$ does not depend on the connection
$A$, since the cohomology classes $\eta_E$ and $[B_2 (F_A , F_A) ]$
are $A$-independent. 
Furthermore, because homotopic sections $\phi :M \rightarrow 
E$ induce the same map $\phi^\ast : H^\ast (E) \rightarrow H^\ast (M)$ on
the cohomology \cite{B-T}, by Stokes theorem $T_{[\phi]}$ only depends
on the homotopy class of  $\phi$.  

$\ $
  
To end this subsection we state two results that, to some extent, clarify the
meaning of the first and the third vortex equations. The first
proposition is well known \cite{MiR}. As for the second proposition,
we relegate its proof to appendix B, since it is a bit long and,
moreover, is just a mild extension of well known calculations
\cite[p. 9]{Ko}.

\begin{prop}
Let $A \in \AAA$ be any connection and let $\phi$ be a section of $E$. Then
$\db^A \phi  = 0$ if and only if  $\phi$ is holomorphic as a map $(M, J_M)
\rightarrow (E, J(A))$.
\label{p2.1}
\end{prop}

\begin{prop}
The condition $F_A^{0,2} = 0$ implies that the almost-complex
structure $J(A)$ on $E$ is integrable. The converse is also true if at
least one point in $F$ has a discrete isotropy group (contained in $G$).
\label{p2.2}
\end{prop}

\subsection{Complex gauge transformations}

Here we recall the notions of complexified Lie group,
complexified action, and complex gauge transformation.
To any compact Lie group $G$ one can associate a complex analytic Lie
group $G_{\CC}$, called the complexification of $G$. The Lie algebra
of $G_{\CC}$ can be identified with the complexification
$\g_{\CC} = \g \oplus i\g $ of the Lie algebra of $G$. Both $G$ and
$\g$ can be naturally embedded into $G_\CC $ and $\g_\CC$,
respectively, as fixed points of natural involutions --- called
conjugations --- in these spaces \cite{B-D}. Furthermore, when the group $G$
acts holomorphically on a compact K\"ahler manifold $F$, this action
can be canonically extended to a holomorphic action of $G_\CC$ on $F$
\cite{G-S}. At a Lie algebra level this extension is defined by
\begin{gather}
(u +iv)^\flat \ =\ u^\flat + J_F \ v^\flat \ \ \qquad \forall \; u,v \in \g \ ,
\label{3.1.0}
\end{gather}
where we denote by $v^\flat$ the  vector field on $F$ defined by the flow 
$t \mapsto \rho_{\exp (t v )}$, and $J_F$ is the complex structure on $F$. 

This extension of the action on $F$ allows us to define complex gauge
transformations on the 
bundle $E=P\times_G F$, which extend to $G_\CC $ the original
$G$-gauge transformations. A complex gauge transformation $g$ is a
section of the bundle $P\times_{{\rm Ad}_G} G_\CC $ over $M$. The set
of these sections forms a group, denoted by ${\mathcal G}_\CC$. Each
$g \in {\mathcal G}_\CC $ determines an automorphism of $E$ by the
formula 
\begin{gather}
[p,q] \ \mapsto\  [p, \rho_{g_p} (q)]\ ,
\label{3.1.1}
\end{gather}
where $\rho$ is the extended $G_\CC$-action and $g_p$ is the only
element of $G_\CC$ such that $g \circ \pi_P (p) = [p, g_p]$. 
If we compose a
section $\phi \in \Gamma (E)$ with this automorphism of $E$ we
get another section, which we denote by $g(\phi )$. Complex gauge
transformations can also be made to act on the space $\AAA$ of
connections on $P$, in such a way as to extend the action of the
original $G$-gauge transformations. This extension is defined by the
formula
\begin{gather}
g(A) \ =\ {\rm Ad}_g \circ A \ -\ \pi_P^\ast (g^{-1} \db g + \bar{g}^{-1}
\partial \bar{g})\ . 
\label{3.1.2}
\end{gather} 
An important fact about these complex gauge transformations is that
both the first and the third vortex equations are invariant by them,
whereas the second equation is invariant by real gauge
transformations only.

For later convenience we also record here the following standard definition.
\vspace{.3cm}
\begin{defn} 
A divisor $D$ on $M$ is a locally finite formal linear combination
\[
D \ =\ \sum_i \ a^i \cdot Z_i \ , \qquad a^i \in {\mathbb Z} \ ,
\]
of irreducible analytic hypersurfaces $Z_i$ of $M$. The divisor $D$ is
called effective if $a_i \geq 0$ for all $i$. The support of $D$, written
${\rm supp}\, D$, is the subset of $M$ formed by the union of the
hypersurfaces $Z_i$ with non-zero coefficient $a^i$.
\end{defn}

\subsection{Torus principal bundles}

In this last subsection we will introduce some notation and recall
some standard results about $T^n$-principal bundles used in the rest
of the paper. 

Let $P \rightarrow M$ be any principal $T^n$-bundle, let $\hat{\rho}$ be the
natural action of $T^n$ on $\CC^n$, denote by $\hat{\rho}_j$ the restriction
of $\hat{\rho}$ to the $j$-th factor $\CC$ in $\CC^n$, and let $\hat{L}_j = P
\times_{ \hat{\rho}_j} \CC $ be the associated line bundle.
We begin by stating a standard result, whose proof we omit.
\begin{prop}
Given any $n$ classes $\alpha_j \in H^2 (M; {\mathbb Z})$ there is
exactly one principal $T^n$-bundle $P \rightarrow M$, up to
isomorphism, such that $\alpha_j$ coincides with the first Chern class
$c_1 (\hat{L}_j)$.
\label{p5.1}
\end{prop}

\vspace{.2cm}

This proposition shows that the correspondence $P \mapsto
\alpha (P)$, with $\alpha_j (P) = c_1 (\hat{L}_j )$, defines a bijection
between the set of principal $T^n$-bundles over $M$ (up to
isomorphism), and the $n$-fold cartesian product of $H^2 (M;{\mathbb
  Z})$. Now identify the Lie algebra $\ttt^n$ with $\RR^n$ in such a
way that the exponential map $\ttt^n \rightarrow \RR^n$ is
\begin{gather}
\exp (w_1, \ldots , w_n ) \ =\ (e^{2\pi i w_1}, \ldots ,e^{2\pi i
  w_n})\ , \qquad w_k \in \RR .
\label{4.1.1}
\end{gather}
With this identification, for any principal $T^n$-bundle $P \rightarrow
M$ we define
\begin{gather}
\deg P \ =\ - \int_M \Lambda F_A \ \: \omega_M^{[m]}\ , 
\label{4.1.2}
\end{gather}
where $A$ is any connection on $P$. This constant does not depend on
$A$. In fact, having in mind the 
above identification of $\ttt^n$ with $\RR^n$, it is clear that   
$2\pi i (F_A)_j$ coincides
with the curvature on the base of the connection on $\hat{L}_j$ induced by
$A$. In particular $c_1 (\hat{L}_j) = - [(F_A)_j]$, and so it follows from
(\ref{2.111}) and proposition \ref{p5.1} that
\begin{gather}
\deg{P}\ =\ -\int_M F_A \wedge \ast \omega_M  \ =\ \int_M \alpha (P)
\wedge \omega_M^{[m-1]} \qquad \in \ \RR^n \ .
\label{4.1.4}
\end{gather}
We also define the constant 
\begin{gather}
c(a,P,M)\ :=\ (a^2\, {\rm Vol}\,M)^{-1} \deg P \ ,
\label{4.1.5}
\end{gather}
which will appear often in the subsequent sections.

$\ $

Finally, to end this subsection, we will state a lemma necessary for
section 6. Let $\beta : T^d \rightarrow T^n$ be any homomorphism of
tori. These have the general form
\[
\beta (g_1 , \ldots , g_d) \ =\ (\ldots , \Pi_{1\le l \le d} \;
(g_l)^{\beta_{al}}  , \ldots )_{1\le a \le n} \ , \qquad {\rm with}\ \ 
\beta_{al} \in {\mathbb Z} \ .
\]
Given a principal $T^d$-bundle $P \rightarrow M$, the associated bundle
$P' = P \times_{\beta} T^n$ is in a natural way a principal
$T^n$-bundle over $M$. Then the following naturality property is easy
to check.
\begin{lem}
The classes in $H^2 (M ; {\mathbb Z})$ associated with $P'$ are
$\alpha_a (P') = \sum_{l=1}^d \beta_{al}\; \alpha_l (P)$ for all $a=1,
\ldots , n$.
\label{p5.2}
\end{lem}

\section{A simpler case: $\CC^n$ with $T^n$-action}

\subsection{Results}

In this section we give the space of solutions and energy spectrum of
the vortex equations in the case $F=\CC^n$ and 
$G=T^n$. The results are contained in theorems \ref{ts1} and
\ref{ts2}. 

$\ $

One starts with the action of $T^n$ on $\CC^n$ given by
\begin{gather}
\rho_{( g_1 , \ldots , g_n )}\; (z_1 , \ldots , z_n ) \ =\ ( \; \cdots ,
z_k \; \Pi_j \: (g_j)^{C_{kj}} , \cdots \; )_{1\le k \le n} \ ,
\label{s.1}
\end{gather}
where the matrix $C$ belongs to $SL(n; {\mathbb Z})$. It is an
effective hamiltonian action. Identifying $\ttt^n \simeq \RR^n$ in the
usual way (\ref{4.1.1}), the general form of a moment map $\mu :
\CC^n \rightarrow \RR^n$ for this action is 
\begin{gather}
\mu (z_1 , \ldots , z_n ) \ =\ - \, \pi\: \left( \cdots ,\ \sum_j C_{jk}\,
|z_j|^2 , \ \cdots \right)_{1\le k \le n} \; + \  t \ ,
\label{s.2}
\end{gather}
where $t$ is any constant in $\RR^n$.

Now consider the associated vector bundle $E = P \times_{\rho} \CC^n
$. Denoting by $\rho_j$ the restriction of the action $\rho$ to the
$j$-th component $\CC$ of $\CC^n$, we have that
\[
E \ =\ L_1 \oplus \cdots \oplus L_n  \ ,
\]
where $L_j = P \times_{\rho_j} \CC $ is the associated line
bundle. Notice that the natural hermitian products on $\CC^n$ and $\CC$
induce hermitian metrics on the bundles $E$ and $L_j$, because the
actions $\rho$ and $\rho_j$ are unitary. We denote by $h$ and $h_j$
these hermitian metrics.

Finally, an integrable connection $A \in {\mathcal A}^{1,1} (P)$
induces a metric-compatible integrable connection $\nabla$
(resp. $\nabla_j$) on the vector bundle $E$ (resp. $L_j$). In
turn, this integrable connection defines a unique holomorphic
structure on $E$ (resp. $L_j$) such that $\nabla$ (resp. $\nabla_j$)
is the hermitian connection of this bundle \cite{Ko}. The bundles $E$
and $L_j$ equipped with these holomorphic structures will be denoted
by $E^A$ and $L_j^A$. Notice that, also as holomorphic hermitian
bundles, 
\begin{gather}
E^A \ =\ L_1^A \oplus \cdots \oplus L_n^A  \ .
\label{s.3}
\end{gather}
Recalling the constants $c(P,M,a) \in \RR^n$ and $\alpha (P) \in H^2
(M; {\mathbb Z})^n$ defined in section 2.4, we have the following
results.
\begin{thm}
In the setting described above, the vortex equations (\ref{2.0}) have
solutions only if the constant $c(P,M,a)$ is in $\mu
(\CC^n)$. When this constant lies in the interior of $\mu (\CC^n)$,
the set of solutions can be described as follows. For each $j = 1 ,
\ldots , n$ pick an effective divisor $D_j = \sum_i a^i_j \cdot Z_i$
on $M$ representing the homology class Poincar\'e dual to $\sum_k
C_{jk} \: \alpha_k (P)$. Then there is a solution $(A , \phi )$ of
(\ref{2.0}) such that $D_j$ is the divisor of the zero set of $\phi_j$
(the $j$-th component of $\phi$ under the decomposition (\ref{s.3}))
regarded as a holomorphic section of $L_j^A$. This solution is unique
up to gauge transformations, and all solutions of (\ref{2.0}) are
obtained in this way. 
\label{ts1}
\end{thm}

\begin{thm}
The topological energy (\ref{2.6}) of any solution of the vortex
equations is
\begin{gather}
T \ =\  \int_M \  \sum_k \  \left( \; t_k \; \alpha_k (P) \wedge
\omega_M^{[m-1]}\ - \ \frac{1}{a^2}\: \alpha_k (P) \wedge \alpha_k (P)
\wedge \omega_M^{[m-2]} \; \right) \ ,
\end{gather}
where $t \in \RR^n$ is the arbitrary constant in the moment map (\ref{s.2}).
\label{ts2}
\end{thm}

In theorem \ref{ts1} it is of course implicit that, if it is
impossible to find a suitable set of divisors $D_j$, then the set of
vortex solutions is empty. Notice as well that the statement of these
results is especially simple when $M$ is a 
Riemann surface, due to the isomorphism $H^2 (M ; {\mathbb Z}) \simeq
{\mathbb Z}$. In fact, in this case it is apparent from theorem
\ref{ts1} that the moduli space of vortex solutions can be identified
with the product of symmetric powers $S^{N_1}M \times \cdots \times
S^{N_n}M$, where $N_j$ is the integer $\sum_k C_{jk}\, \alpha_k (P)$. If
any of these integers is negative, then the moduli space is empty. For
$M$ a Riemann surface the topological energy also reduces to $T =  t
\cdot \alpha (P)$. Another interesting fact regarding the
topological energy is that, unlike the $\CC {\mathbb P}^n$ case of
theorem \ref{t6.1}, here the energy is completely determined by the
bundle $P$; it does not depend on the particular solution chosen. This
difference between the $\CC^n$ and the $\CC {\mathbb P}^n$ cases is
analogous to the fact that the degree of a line bundle completely
determines the number of zeros of a holomorphic section, but not of a
meromorphic section.

$\ $

The key ingredient to prove theorem \ref{ts1} is the following
proposition, which follows quite straightforwardly from the
``stability'' criterion of \cite{Ba}.
\begin{prop}
Assume that $c(P,M,a)$ lies in the interior of $\mu (\CC^n)$, and let
$(A, \phi) \in {\mathcal A}^{1,1} (P) \times \Gamma (E)$ be any pair
such that $\db^A \phi = 0$ and $\phi_j$ is not identically zero for
any $j$. Then there exists a complex gauge transformation $g : M
\rightarrow (\CC^\ast )^n $, unique up to multiplication by real gauge
transformations, such that the pair $(g(A), g(\phi))$ is a solution of
the vortex equations. 
\label{ps3}
\end{prop}

When the target of the $\si$-model is a compact manifold, instead of
$\CC^n$, things get rather more complicated. In the next section we
will try to find results analogous to proposition \ref{ps3} in the
compact setting. For this we will use results of \cite{A} and a more
general ``stability'' criterion of \cite{MiR}.

\subsection{Proofs}

The proofs below may be regarded as a warm up to the calculations of
sections 4  and 5. Nevertheless, in order to avoid repetition, we
will occasionally  invoke results from those sections.

\noindent
{\bf Proof of proposition \ref{ps3}.}
The proof is based on the Hitchin-Kobayashi correspondence of
\cite{MiR, Ba}. In the case $G= T^n$ and $X=\CC^n$ this correspondence reduces
to the following statement.

{\it Given a simple pair $(A, \phi) \in {\mathcal A}^{1,1} \times
  \Gamma (E)$, there exists a complex gauge transformation that takes
  this pair to a solution of (\ref{2.0b}) iff
\begin{gather}
v\; \cdot \; (\deg P + a^2 ({\rm Vol} M) \, t ) \ < 0  \qquad
{\rm for \ all}\ v \in \RR^n \setminus \{ 0 \} \ {\rm such\ that\ }
(Cv)_j \ge 0 \ .
\label{s.4}
\end{gather}
When it exists, this transformation is unique up to composition with
real gauge transformations.}
 
Now notice that condition (\ref{s.4}) is equivalent to 
\[
v \; \cdot \; (C^{-1})^T (c(P,M,a) + t ) \ < \ 0  \qquad
{\rm for \ all}\ v \in \RR^n \setminus \{ 0 \} \ {\rm such\ that\ \; }
v_j \ge 0 \ ,
\]
or in other words, to the condition that the constant $(C^{-1})^T
(c(P,M,a) + t ) $ lies in the set 
\[
\{ x \in \RR^n : \ x_j < 0 \ \ \forall j     \} \ .
\]
But this is the same as demanding that $c(P,M,a)$
should lie in the interior of $\mu (\CC^n)$, and this is satisfied by
assumption.  Thus in order to prove the proposition it is enough to
show that $(A, \phi)$ is a simple pair. This can be done just as in
the proof of Proposition \ref{p4.7}.  \hfill \qed

\vspace{.4cm}

\noindent
{\bf Proof of theorem \ref{ts1}.}
The first statement of the theorem can be proved by integrating the
second vortex equation over $M$, just as in the proof of Theorem
\ref{t4.1}.

For the rest of the theorem, assume that $c(P,M,a)$ lies in the
interior of $\mu (\CC^n)$, and consider the divisors $D_j$ described
in the theorem. Since
\[
L_j \ =\ \bigoplus_{1 \le j \le n} (\hat{L}_j)^{C_{kj}} \ ,
\]
where the line bundles $\hat{L}_j$ were defined in section 2.4, we
have that
\[
{\rm PD}(D_j) \ =\ \sum_k C_{jk} \: \alpha_k (P) \ =\ \sum_k C_{jk} \:
c_1 (\hat{L}_k ) \ =\ c_1 (L_j) \ .
\]
So it follows from well known results that the divisor $D_j$
determines a holomorphic structure on $L_j$ together with a
non-zero holomorphic section $\phi_j$ of this bundle such that $D_j$
is the zero set divisor of $\phi_j$ \cite{G-H}. Denoting by $\nabla_j$
the hermitian connection of $(L_j , h_j)$ equipped with this
holomorphic structure, by construction we have that $\nabla_j^{0,1}
\phi_j =0 $. Now, just as in the first part of the proof of Lemma
\ref{l6.8}, there exists a connection $A$ on $P$ such that $\nabla_j$
is the connection on $L_j = P \times_{\rho_j} \CC $ induced by $A$. So
using the decomposition (\ref{s.3}) to define
\[
\phi \ := \ (\phi_1 , \ldots , \phi_n ) \quad \in \ \Gamma (E) \ ,
\] 
we have that $\db^A \phi = (\ldots , \nabla_j^{0,1} \phi_j , \ldots )
=0$. The existence part of Theorem \ref{ts1} then follows from
Proposition \ref{ps3} together with the fact that complex gauge
transformations do not change the zero set divisor of a section (see
Lemma \ref{l6.9}).

As for the unicity of the solutions, suppose that $(A_1 ,  \phi_1)$ and
$(A_2 , \phi_2)$ are two solutions of the vortex equations such that the
zero set divisors of $\phi_1$ and $\phi_2$ are well defined and
equal. Then Lemma \ref{l6.9} tells us that the components $(\phi_1)_j$
and $(\phi_2)_j$, and therefore also $\phi_1$ and $\phi_2$, are
complex-gauge equivalent. The unicity statement of Proposition
\ref{ps3} then guarantees that $\phi_1$ and $\phi_2$ are real-gauge
equivalent, as required.

Finally, to recognize that all solutions of the vortex equations are
of the kind described in Theorem \ref{ts1}, it is enough to show that
if $(A, \phi)$ is a solution then the zero set divisor of $\phi_j \in
\Gamma (L_j^A )$ is well defined, that is $\phi_j$ is not the zero
section for any $j$. But if this were not true, the image $\mu \circ
\phi (M)$ would be contained in one of the boundary faces of $\mu
(\CC^n)$, and then integrating the second vortex equation over $M$ one
would obtain a contradiction with the fact that $c(P,M,a)$ does not
belong to this face (see the analogous result in Theorem \ref{t4.1}).
\hfill  \qed

\vspace{0.4cm}

\noindent
{\bf Proof of theorem \ref{ts2}.}
Let $(A, \phi)$ be any solution of the vortex equations. We start by
considering the 1-form on $P\times \CC^n$
\[
\nu \ = \ \frac{i}{2} \: \sum_k z^k \: \dd \bar{z}^k \ - \ (\mu -
t)\cdot A \ .
\]
It is not difficult to check that this form annihilates vectors in the
kernel of the quotient map $\chi : P\times \CC^n \rightarrow E$, and
that it is invariant under the right action $(p,\: v)\cdot g = (p\cdot g
, \: \rho_{g^{-1}} v)$  of $T^n$ on $P\times \CC^n$. Thus $\nu$
descends to a form on $E$, that is $\nu = \chi^\ast \nu_E$ for some
complex 1-form $\nu_E$ on $E$. On the other hand, using
(\ref{2.9}), we have that 
\[
\dd \nu \ =\ \frac{i}{2} \sum_k \dd z^k \wedge \dd \bar{z}^k \ -\ \dd
(\mu \cdot A ) \ + \ t \cdot \dd A \ = \ \eta (A) \ +\ t\, \cdot\,
\pi^\ast_P F_A \ , 
\]
where $F_A$ is the curvature form on the base $M$ of the connection
$A$. But by the commutativity of the diagram
\begin{equation*}
\begin{CD}
P \times \CC^n        @>{\chi}>>      E             \\ 
@VVV                                  @VV{\pi_E}V   \\
P                     @>>{\pi_P}>     M
\end{CD}
\end{equation*}
we get that on $P\times \CC^n$
\[
\pi_P^\ast \: F_A \ = \ \chi^\ast \, \pi_E^\ast \: F_A \ ,
\]
and so
\[
\chi^\ast \: \dd\, \nu_E \ =\ \dd \,\nu \ =\  \chi^\ast (\: \eta_E (A) +
t\cdot \pi^\ast_E F_A \: ) \ .
\]
This implies that
\[
\eta_E (A) \ =\ \dd\, \nu_E \ -\ t \,\cdot\, \pi_E^\ast\, F_A \ ,
\]
and using Stokes theorem we get that
\[
\int_M \phi^\ast \,\eta_E (A) \: \wedge \: \omega_M^{[m-1]} \ =\ - \: t \:
\cdot \: \int_M F_A \wedge \omega_M^{[m-1]} \ .
\]
The formula for the topological energy then follows from the
definition (\ref{2.6}) and the identification $\alpha_j (P) = -
[(F_A)_j]$ in $H^2 (M)$ given in section 2.4. 
\hfill   \qed

\section{The ${\bf 2^{\rm nd}}$ vortex equation as an imaginary-gauge fixing
  condition}

\subsection{Main results}

As was mentioned in section 2.3, an important
fact about the complex gauge transformations is that
both the first and the third vortex equations are invariant by
them, whereas the second equation is not. Hence, given a pair
$(A,\phi)$ that solves (\ref{2.0a}) and (\ref{2.0c}), it makes sense
to ask whether 
there is a complex gauge transformation $g$ such that $(g(A), g(\phi
))$ solves (\ref{2.0b}), and therefore all the vortex equations
\cite{MiR}. The 
ideal answer would be that such a transformation always exists and is
unique up to real gauge transformations. This would mean that equation
(\ref{2.0b}) acts as a sort of imaginary-gauge fixing condition, and that the
set of solutions of (\ref{2.0}) up to real gauge transformations is
the same as the set of solutions of (\ref{2.0a}) and
(\ref{2.0c}) up to complex gauge
transformations.

The purpose of this section is to study this problem when the gauge
group is $T^n$ and the target $F$ is compact.
 The basic results obtained are expressed in
theorems \ref{t4.1} and \ref{t4.2}. As a kind of corollary we find
that although 
the ideal answer stated above is not in general true, it comes very
close to being completely true when the target manifold $F$ is a
``simple'' one --- for example when $F$ is toric (see corollary \ref{c4.3} and
the following remark). This will eventually allow us to compute the
moduli space of solutions when $F= \CC {\mathbb P}^n $ (section 5), and to find
a big set of non-trivial solutions for more general toric $F$'s (section 7).

$\ $

In order to state the basic results of this section we first need to
establish some notation. The complexified torus is $T^n_\CC \simeq
(\CC^\ast )^n$, and its Lie algebra is identified with $\ttt^n
\oplus i\ttt^n \simeq \CC^n $ in such a way that the
exponential map is 
\begin{gather}
\exp (w_1, \ldots , w_n ) \ =\ (e^{2\pi i w_1}, \ldots ,e^{2\pi i
  w_n})\ , \qquad w_k \in \CC .
\label{5.1.1}
\end{gather}
The inner product on $\ttt^n \simeq \RR^n$ is just the euclidean
one. For any point $p$ in $F$ we call $\OO_p$ and $\OO_p^\CC $ its
$T^n$-orbit  and $T^n_\CC$-orbit, respectively; similarly, the real and
complex isotropy groups of $p$ are denoted by $G_p$ and $G_p^\CC
$. 

Also a word about polytopes. By the convexity theorem, if
$\mu : F \rightarrow \ttt^n \simeq \RR^n $ is a moment map for a torus
action on $F$, which is assumed compact, its image $\mu (F)$ is a
convex polytope in $\RR^n$ 
(see for example \cite{A} or \cite{MD-S}). As
a set, $\mu (F)$ is the disjoint union of its $k$-dimensional open
faces, or $k$-cells, for $k = 0, \ldots, \dim \mu (F)$. Thus for
example $\mu (F)$ has only one open face of maximal dimension, and the
$0$-dimensional open faces are the vertices of $\mu (F)$. We are now
ready to state the main results of this section.

\begin{thm}
A necessary condition for the equation $\Lambda F_A + a^2\ \mu\circ\phi
=0$ to have a solution is that the constant $c(a, P, M)$
(c.f. (\ref{4.1.5})) lies in $ \mu (F)$. If this is 
satisfied, let $\si_c$ be the only open face of the polytope $\mu (F)$
that contains this point, and let $\sib_c$ denote its closure. Then for
any $(A,\phi )$ solution of (\ref{2.0b}), the image $\mu \circ \phi (M)$ is
contained in $\sib_c$ and is not entirely contained in any of the closed faces
of  $\sib_c \setminus \si_c $.
\label{t4.1}
\end{thm}

\begin{thm}
Let $(A,\phi ) \in \AAA \times \Gamma (E)$ be a pair such that, for all $x$
in some open dense subset of $M$, the conditions
\begin{itemize}
\item[(i)] $\OO_{\phi (x)}^\CC  \ \cap \ \mu^{-1}(c) \ \ne \emptyset \ $;
\item[(ii)] $G_{\phi (x)}$ has dimension $n - \dim{\si_c} \ $;
\end {itemize}
are satisfied. Then there exists a complex gauge transformation $g : M
\rightarrow T_\CC^n$ that takes $(A, \phi )$ to a solution of
(\ref{2.0b}). This 
transformation is unique up to multiplication by transformations whose
imaginary part is a constant in $\exp (i
\si_c^{\perp})$.
\label{t4.2}
\end{thm}

\begin{rem}
Here we will only prove this theorem in the generic case where the
constant $c(P,M,a)$ lies in the interior of $\mu (F)$, i.e. when $\dim
\si_c = n$. This is the only case needed in the subsequent sections.
The proof is based on a very general criterion of \cite{MiR}. For a
hint of the proof in the general case see the remark in section 4.3.
\end{rem}

\begin{cor}
Assume that the orbit $\OO_p^\CC$ of any point $p \in \mu^{-1}(\si_c ) $
satisfies $\mu (\OO_p^\CC ) = \si_c$. Then, given any pair $(A,\phi )\in
\AAA^{1,1} \times \Gamma (E) $ such that $\db^A \phi = 0$, there exists a
complex gauge transformation that takes $(A, \phi )$ to a solution of
(\ref{2.0b}) 
if and only if the image $\mu \circ \phi (M)$ is contained in $\sib_c$ but
not in any of the closed faces of $\sib_c \setminus \si_c$. Furthermore,
when it exists, this transformation is unique up to multiplication by
transformations whose imaginary part is a constant in $\exp (i
\si_c^{\perp})$.
\label{c4.3}
\end{cor}

\begin{rem}
The condition of this corollary, namely 
$\mu (\OO_p^\CC ) = \si_c $ for any $p \in \mu^{-1}(\si_c ) $, is
very restrictive. It is satisfied, however, when the action is
effective and $\dim_{\CC}F =n= \dim_{\RR}T^n$. In this case $F$  
becomes a compact K\"ahler toric manifold, and it is well known that for such
manifolds there is a one-to-one correspondence between open faces of $\mu
(F)$ and $T^n_{\CC}$-orbits in $F$, which is given by $\si \mapsto
\mu^{-1}(\si )$.
This is valid for all K\"ahler toric manifolds, not just the canonical
ones described in section 7.1. 
\end{rem}

\subsection{Proof of theorem \ref{t4.1} and corollary \ref{c4.3}}

We begin this subsection with the proof of theorem \ref{t4.1}. After
establishing two auxiliary lemmas, we end it with the proof of
corollary \ref{c4.3}.

\vspace{0.5cm}

\noindent
${\bf Proof\ of\ theorem\ \ref{t4.1}.}$ 
Let $(A,\phi )$ be a solution of (\ref{2.0b}). Integrating this
equation over $M$ and using (\ref{4.1.2}) and (\ref{4.1.5}) one has that
\begin{gather}
\int_M (\mu \circ \phi - c)\ \omega_M^{[m]} \ \ =\ \ 0 \qquad \in \ \
\RR^n .  
\label{5.2.1}
\end{gather}
If $c \not\in \mu (F)$, from the convexity of $\mu (F)$ it is clear that
for all $v\in \mu (F)$ the vectors $v-c$ will lie in the same open
half-space of $\RR^n$. In particular the same thing happens with the
vectors $\mu \circ \phi (x) - c $ for all $x\in M$, and thus it is
impossible for (\ref{5.2.1}) to hold --- a contradiction.

Now suppose that $c$ lies in some open face $\si $ of $\mu (F)$. If  $\si
$ is $n$-dimensional, it is obvious that $\mu \circ \phi (M) \subseteq
\sib = \mu (F)$. If the dimension of $\si $ is $k < n$, let $A_1 , \ldots
, A_{n-k}$ be the closed $(n-1)$-dimensional faces of $\mu (F)$ whose
intersection is 
$\sib $, and let $n_j$ be an outward normal vector to $A_j$. Then, from the
convexity of $\mu (F)$, one has that $n_j \cdot (v-c) \leq 0 $ for all $v \in
\mu (F)$, and the equality holds iff $v \in A_j $. But (\ref{5.2.1})
implies that 
\[
\int_{x\in M} n_j \cdot (\mu \circ \phi (x) - c) \ \; \omega_M^{[m]} \
= \ 0 \ , 
\]
and so we conclude that $\mu \circ \phi (x) \in A_j$ for all $x \in
M$. Since this is true for all $j$ we actually have that $\mu \circ
\phi (M) \subseteq \sib $, as required.

On the other hand, let $B$ be any closed face of $\sib \setminus \si$ ---
which is also a $(k-1)$-dimensional closed face of $\mu (F)$ --- and let
$u$ be a vector normal to $B$, parallel to $\sib $, and pointing outward 
of $\mu (F)$. Then, because $c \in \si $ and $\sib $ is convex, one has
that $u \cdot (v-c) > 0$ for all $v\in B$. In particular it is
impossible that $\mu \circ \phi (M) \subseteq B$, otherwise one would
have that
\[
\int_{x \in M}  u \cdot  (\mu \circ \phi (x) -c) \ \; \omega^{[m]}_M
\ \ > \ \ 0 \ , 
\]
which contradicts (\ref{5.2.1}).  \hfill  \qed

\begin{lem} 
Let $\si $ be any open face of the polytope $\mu (F)$, and denote by $\sib
$ its closure. Then
\begin{itemize}
\item[(i)] $\mu^{-1} (\sib ) $ is a connected complex submanifold of $F$;
\item[(ii)] $\mu^{-1} (\si )$ is invariant under the $T_{\CC}^n$-action.
\end{itemize}
\label{l4.4}
\end{lem}
This lemma is a well known result. Statement (i) follows rather
straightforwardly from lemmas $5.53$ and $5.54$ of \cite{MD-S} and
their proof; statement (ii) follows from Theorem 2 of \cite{A}.

\begin{lem}
Let $\phi : M \rightarrow E$ be a section of $E$ such that $\db^A \phi =
0$ for some connection $A \in \AAA^{1,1} (P)$. Then for any open face $\si $ of
$ \mu (F)$ the inverse image $(\mu \circ \phi )^{-1} (\sib )$ is an
analytic subvariety of $M$.
\label{l4.5}
\end{lem}

\begin{prooff}
To avoid any confusion, in this proof we will use different symbols for
the moment map $\mu : F \rightarrow \RR^n$ and its lift $\tilde{\mu}: E
\rightarrow \RR^n$.
To start with, notice that $\sib $ is a disjoint union of open faces of  
$\mu (F)$, possibly with different dimensions, and so it follows from
lemma \ref{l4.4} that $\mu^{-1} (\sib )$ is a complex submanifold of $F$ which
is invariant by the $T^n$-action. It is not difficult to check that this
implies that  $E' = P \times_{T^n} \mu^{-1}(\sib )$ is a complex
submanifold of $E=P\times_{T^n} F$, where the complex structure on $E$ is
$J(A)$. Furthermore, from the definition $\tilde{\mu}\circ \chi  (p,q) =  
\mu (q)$ (see section 2.1), we also have that $E' = \tilde{\mu}^{-1} (\sib
)$.

On the other hand, by proposition \ref{p2.1}, the section $\phi $ is a
holomorphic map from $M$ to $(E, J(A))$. This map is proper because $M$ is
compact, and since a section is always an immersion, we conclude that $\phi
(M)$ is actually a complex submanifold of $E$, and $\phi : M \rightarrow
\phi (M)$ is a biholomorphism. It is then clear that $E' \cap \phi (M)$,
being an intersection of complex submanifolds, is an analytic subvariety
of $\phi  (M)$. Hence $\phi^{-1} (\tilde{\mu}(\sib )) = \phi^{-1} (E' \cap
\phi (M)) $ is an analytic subvariety of $M$.
\end{prooff}

\vspace{.2cm}

\noindent
{\bf Proof of corollary \ref{c4.3}.} 
We denote by $(\ast)$ the condition ``$\mu \circ \phi (M)$ is
contained in $\sib_c$ but not entirely in any of the closed faces of $\sib_c
\setminus \si_c$''. The proof of necessity is fast, due to theorem
\ref{t4.1}. In fact, it follows from part (ii) of lemma \ref{l4.4}
that, if a section $\phi \in \Gamma (E)$ satisfies $(\ast)$, so does
its entire complex gauge equivalence class. (We are using that the
closed faces of $\sib_c$ are themselves a union of open faces of
$\mu (F)$, and so their inverse image by $\mu$ is also
$T^n_\CC$-invariant.) The necessity of $(\ast)$ is then a direct
consequence of theorem \ref{t4.1}.

To prove the sufficiency and uniqueness statements we will use theorem
\ref{t4.2}. If $\db^A \phi =0$ and $(\ast)$ is satisfied, by lemma
\ref{l4.5} it is true that for each closed face $B$ of $\sib_c \setminus
\si_c$, the inverse image $(\mu \circ \phi)^{-1} (B)$ is an analytic
subvariety of $M$ which is not the entire $M$; in particular this set has
zero measure in $M$. Since this is true for all the faces of $\sib_c
\setminus \si_c $, we conclude that $(\mu \circ \phi)^{-1} (\si_c)$ is
open and dense in $M$. Now if $x \in (\mu \circ \phi)^{-1} (\si_c)$,
that is $\phi (x) \in \mu^{-1} (\si_c)$, by assumption $\mu
(\OO^{\CC}_{\phi (x)}) = \si_c $. Hence on the one hand, since $c(P,M,a)
\in \si_c$, this implies that condition (i) of theorem \ref{t4.2} is
satisfied; on the other hand, using lemma \ref{l4.6}, this implies
that condition (ii) of theorem \ref{t4.2} is satisfied as
well. Applying this theorem we obtain the sufficiency and
uniqueness parts.      \hfill          \qed

\subsection{Proof of theorem \ref{t4.2}}

We first derive an auxiliary lemma and then prove theorem \ref{t4.2}
in the case where $c(P,M,a)$ lies in the interior of $\mu (F)$. At the
end of the subsection we make a remark about the proof in the general
case.

The lemma is the following. Given $p \in F$, let $\si_p$ denote the
only open face of $\mu (F)$ that contains the point $\mu (p)$. Then
by lemma \ref{l4.4} and theorem $2$ of \cite{A}, the image $\mu 
(\OO^{\CC}_p )$ is a convex open polytope contained in $\si_p$.

\begin{lem}
Given $p \in F$, the Lie algebra of the isotropy subgroup $G_p \subseteq
T^n$ is the subspace of $\ttt^n$ formed by the vectors orthogonal to $\mu
(\OO^{\CC}_p)$. In particular ${\rm Lie }\, G_p$ contains the subspace
$\si_p^\perp $.
\label{l4.6}
\end{lem}

\begin{prooff}
Theorem 2 of \cite{A} guarantees that the restriction of $\mu $ to
$\OO^{\CC}_p $ induces a homeomorphism $\OO^{\CC}_p  / T^n  \rightarrow
\mu (\OO^{\CC}_p ) $. Since the dimension of the isotropy subgroup $
G^{\CC}_p \subseteq T_{\CC}^n $ is twice the dimension of $G_p$, we    
conclude that $\OO^{\CC}_p  / T^n $, and therefore $\mu (\OO^{\CC}_p )$,
have dimension $n-\dim{G_p}$. On the other hand, for any $v\in
\ttt^n$, property $(i)$ of the definition of a moment map (see section
2.1) implies that
\[
v \ \perp\ {\rm Image}(\dd \mu)_p  \quad  \iff \quad  (v^\flat)_p\ =\ 0
\quad  \iff  \quad  v\ \in \ {\rm Lie }\, G_p \ .
\]
Since
\[
{\rm Image}(\dd\, \mu )_p \quad \supseteq \quad (\dd\, \mu )_p (T_p
\OO^{\CC}_p ) \ =\ T_{\mu (p)} \; \mu (\OO^{\CC}_p ) \   ,
\]
after identifying  $T_{\mu (p)} \ttt^n \simeq \ttt^n$ we obtain that
Lie$\, G_p$ is contained in the subspace of $\ttt^n$ orthogonal to $\mu
(\OO^{\CC}_p )$. Comparing the dimensions, we conclude that Lie$\, G_p$ is  
in fact equal to that subspace.
\end{prooff}

\begin{prop}
Theorem \ref{t4.2} is true when $\dim{\si_c} = n$.
\label{p4.7}
\end{prop}

\begin{prooff}
To prove this proposition we will use the results of \cite{A} and the
Hitchin-Kobayashi correspondence of \cite{MiR}. The latter result is hugely
simplified for abelian $G$, which is the case that matters to us, and can 
be stated in the following form \cite{MiR}.

{\it  Given a simple pair $(A, \phi ) \in \AAA^{1,1}\times \Gamma (E)$,  
there exists a complex gauge transformation that takes this pair to a
solution of (\ref{2.0b}) iff
\begin{gather}  
-v\cdot {\rm deg }P \ +\ a^2 \int_{x\in M} \lambda (\phi (x) , v) \quad
> \ \ \ 0 \quad \ {for\ all\ } v\in \RR^n \ .
\label{5.3.1}
\end{gather}
When it exists, this transformation is unique up to composition with real
gauge transformations. }
  
Thus to prove the lemma we only have to show that the pair $(A ,\phi )$ of
theorem \ref{t4.2} is simple and satisfies (\ref{5.3.1}). The
definition of the function
$\lambda $ under the integral is the following. Let $\eta_t^v : F
\rightarrow F $ be the gradient flow of the function $v\cdot \mu : F  
\rightarrow \RR $, and write $\phi (x) = \chi (p,q)$ (see section 2.1); then
\[
\lambda (\phi (x) , v) \ := \ \lim_{t\rightarrow +\infty} \ v\cdot \mu
(\eta_t^v (q))\ .
\]
The integral of $\lambda$ over $M$ is not in general an easy number to
estimate. However, when the assumptions of theorem \ref{t4.2} hold,
this obstacle 
evaporates, and we will now see how. Take $x\in M$ such that (i) and (ii)
hold. By (i) the constant $c$ is in the open polytope $Q_x := \mu
(\OO^{\CC}_{\phi (x)})$; by (ii) and lemma \ref{l4.6}, this polytope has
dimension $n$. Therefore using lemma $3.1$ of \cite{A}, we have that

\[
\lim_{t\rightarrow +\infty} \ v\cdot \mu (\eta_t^v (q)) \ = \ \sup_{p \in
\OO^{\CC}_{\phi (x)}} v \cdot \mu (p) \ = \ \sup_{u \in Q_x} v\cdot
u \ > \ v\cdot c \ ,
\]
where the strict inequality follows from $Q_x$ being  open and having
dimension $n$ (in particular $v$ cannot be orthogonal to $Q_x$). Since
this holds for all $x$ in an open dense subset of $M$, we conclude that
\[
\int_{x\in M} \lambda (\phi (x) , v) \quad > \ \ \ ({\rm Vol}\ M) \
v\cdot c \qquad {\rm for\ all}\ v\in \RR^n \ ,
\]
which is equivalent to (\ref{5.3.1}).

To prove that the pair $(A, \phi)$ is simple (for the definition of simple
pair see \cite{MiR}), it is enough to show that any infinitesimal gauge
transformation $s : M \rightarrow {\rm Lie}\, T^n_{\CC}=\CC^n $ that leaves
$(A, \phi)$ fixed is necessarily zero. Let $\AAA$ be the space of   
connections on $P$, and as usual identify $T_A \AAA \simeq \Omega^1 (M,
\ttt^n)$. The infinitesimal gauge transformation $s$ produces a tangent
vector in $T_A \AAA$, and by the transformation rules (\ref{3.1.2})
this is given by 
\[
\db s \ +\ \partial \bar{s} \qquad \in\ \Omega^1 (M, \ttt^n) \ .
\]
But if $s$ leaves $A$ fixed, that is $\db s +\partial \bar{s}  = 0$, the
decomposition $\Omega^1 = \Omega^{1,0} \oplus \Omega^{0,1}$ implies that
$\db s = 0$, and since $M$ is compact the function $s$ must be constant.
On the other hand for any $x\in M$ such that (ii) is satisfied we have
that ${\rm Lie}\, G_{\phi (x)} = \{ 0\} $, and so $s(x)$ leaves $\phi (x)
\in E $ fixed iff $s(x) = 0$. By the constancy of $s$ we finally conclude
that $s=0$, and this finishes the proof.
\end{prooff}

\begin{rem}
The proof for the general case $\dim \si_c \leq n$ goes along the
following lines. If $c(P,M,a)$ lies in the boundary of $\mu (F)$, by
assumption (i) and lemma \ref{l4.4} so does the image $\mu \circ \phi
(M)$. Lemma 
\ref{l4.6} then tells us that there is a subtorus of $T^n$ that acts
trivially on $\phi$. The strategy of the proof is to eliminate this
subtorus by formulating the problem in terms of the quotient group and
quotient principal bundle. The isotropy groups $G_{\phi (x)}$ of
assumption (ii) will then have dimension zero, and we will be reduced
to the case $\dim \si_c = n$.
\end{rem}

\section{The vortex solutions for target $\CC {\mathbb P}^n$ }

\subsection{The main result}

We start with the natural action of $T^{n+1}$ on $\CC {\mathbb P}^n$,
given in homogeneous coordinates by 
\[
(g_0 , \ldots ,g_n) \cdot [z_0 , \ldots , z_n] \ =\ [g_0 z_0
  , \ldots , g_n z_n ] \ .
\]
Although this is not an effective action, it induces an effective
hamiltonian action of the quotient group $T^{n+1}/N$, where $N$
denotes the diagonal circle inside $T^{n+1}$. Now, this quotient group
is isomorphic to $T^n$ but, since there is no canonical choice of
isomorphism, there are several different ways of implementing the
$T^{n+1}/N$-action as an action of $T^n$ on $\CC {\mathbb P}^n$. The
general formula for these $T^n$-actions is  
\begin{gather}
\rho_{(g_1 ,\ldots ,g_n)}\, (\,[z_0, \ldots , z_n]\,) \ =\ [z_0 \, ,\: 
 \Pi_j (g_j)^{C_{1j}} z_1 \, ,  \ldots , \: \Pi_j (g_j)^{C_{nj}} z_n 
  \: ] \ ,
\label{6.1}
\end{gather}
where the matrix $C$ is in $SL(n, {\mathbb Z})$. The different choices
 of $C$ correspond to the different possible isomorphisms $T^{n+1}/N
 \simeq T^n$. These actions are all hamiltonian and, using the
identification $(\ttt^n)^\ast \simeq \ttt^n \simeq \RR^n $ 
determined by (\ref{4.1.1}), the general form of a moment map $\mu :
 \CC {\mathbb P}^n \rightarrow \RR^n$ is
\begin{gather}
\mu ([z_0 , \ldots , z_n ])\ =\ \frac{-\: \pi}{\sum_i |z_i|^2 }\; (\ldots ,
\sum_{j\ge 1} C_{jk} |z_j|^2 , \ldots )_{1\leq k \leq n}\ +\ {\rm
  const.} \ .
\label{6.mm}
\end{gather}
We denote by $\Delta$ the image $\mu (\CC {\mathbb P}^n )$ in
$\RR^n$. This is clearly a convex polytope, since it is the image of
the standard $n$-simplex
\[
\{ x \in \RR^n :\ 0\leq x_k \leq 1 \ \ {\rm and}\ \ \Sigma_k \; x_k \leq 1 \}
\]
by the invertible linear transformation $-\pi\; C^T$, possibly composed with a
translation.

$\ $

The aim of this section is to prove theorems \ref{t6.1} and \ref{t6.2},
stated below. They characterize the space of solutions and energy spectrum,
respectively, of the vortex equations for target
$\CC {\mathbb P}^n$ with the $T^n$-action described above. Also theorem
\ref{t6.3}, which appears here as an intermediate step to prove
theorem \ref{t6.1}, may have some independent interest.
Before stating these theorems, however, some notation must be
introduced.

$\ $

Let $B_0 , \ldots , B_{n}$ be the $(n-1)$-dimensional faces of the
polytope $\Delta$. We denote by $\beta_j \in {\mathbb Z}^n$ the unique
primitive normal vector to $B_j$ that points to the exterior of
$\Delta$. For each $j = 0, \ldots , n$ define
\[
F_j \ =\ \mu^{-1} ( B_j ) \ =\ \left\{ [z_0, \ldots , z_n ] \in \CC
{\mathbb P}^n :\ z_j = 0\  \right\} \ ,
\]
which is a $\CC {\mathbb P}^{n-1}$ inside $\CC {\mathbb P}^{n}$.
 Since $F_j$ is a $T^n$-invariant complex submanifold of
$\CC {\mathbb P}^n$, as in the proof of lemma \ref{l4.5} one can
 define the sub-bundles $E_j = P 
\times_{\rho} F_j $ of $E$; these are complex submanifolds of $(E,
J(A))$, where $J(A)$ is the complex structure on $E$ induced by an
 integrable 
connection $A$ on $P$. Recalling also the constants $c(P,M,a)\: \in
 \RR^n$ and $\alpha (P)\: \in H^{2} (M; {\mathbb Z})^n$ defined in
 section 2.4, we have the following results.

\begin{thm}
In the setting described above, the vortex equations (\ref{2.0}) have
solutions only if the constant $c(P,M,a)$ is in $\Delta$. When
this constant lies in the interior of $\Delta$,  the set of solutions
can be described as follows. For 
each $j=0, \ldots , n$ pick an effective divisor $D_j = \sum_i a^i_j
\cdot Z_i$ on $M$ such that 
\begin{itemize}
\item[(i)] the intersection of hypersurfaces ${\rm supp}\, D_0 \cap
  \cdots \cap {\rm supp}\, D_n$ is empty;
\item[(ii)] the Poincar\'e duals (PD) of the fundamental homology cycles
  carried by the divisors $D_j$ satisfy $\alpha (P) = \sum_j \beta_j
  \, {\rm PD}(D_j)$ in $H^2 (M ; {\mathbb Z})^n$.
\end{itemize}
Then there is a solution $(A,\phi )$ of (\ref{2.0}), unique up to
gauge equivalence, such that the intersection multiplicities of the
complex submanifolds $\phi (M)$ and $E_j$ satisfy
\[
{\rm mult}_{\phi (Z_i)} (E_j , \phi (M))\ =\ a^i_j \ .
\]
Furthermore all the solutions of (\ref{2.0}) are obtained in this way. 
\label{t6.1}
\end{thm}

\begin{thm}
Assume that $c(P,M,a)$ lies in the interior of $\Delta$, and let $(A,
\phi)$ be a solution of the vortex equations characterized by divisors
$D_j$, as in the theorem above. Then the topological energy (\ref{2.6}) of
this solution is 
\[
T_{[\phi]} \ =\ e(P,M,\mu , a) \ +\ \frac{\pi}{n+1}\ \sum_{j=0}^{n}\
\int_{M}  {\rm PD}(D_j) \wedge \omega_M^{[m-1]} \ ,
\]
where the constant $e$ does not depend on $(A, \phi)$. Denoting by $b
\in \RR^n$ the barycentre of the polytope $\Delta$, the value of this
constant is
\[
e \ = \ \sum_{k=1}^{n} \ \int_M  \left\{ \; b_k \ \alpha_k (P) \wedge
\omega_M^{[m-1]} \ -\ \frac{1}{a^2} \, \alpha_k (P) \wedge \alpha_k
(P) \wedge \omega^{[m-2]}_M  \; \right\} \ .
\] 
\label{t6.2}
\end{thm}
\vspace{-.6cm}
\begin{rem}
The statement of these results is especially simple when $M$ is a
Riemann surface, for in this case the hypersurfaces $Z_i$ are just points
and, under the isomorphism $H^2 (M; {\mathbb Z}) \simeq {\mathbb Z}$, there
is an identification ${\rm PD}(D_j) \simeq \sum_i a^i_j \ \in {\mathbb
Z}$. In fact, consider the symmetric products $S^{N_j}M$ of the
surface. Each point in $S^{N_j}M$ is an unordered multiplet $(p_1 ,
\ldots , p_{N_j})$ of points in $M$. Now let $\Sigma_{N_0 , \ldots ,
N_n}(M)$ denote the open  dense subset of $S^{N_0}M \times \cdots
\times S^{N_n}M$  obtained by excluding the points that contain a
common $p\in M$ in all the $n+1$ multiplets corresponding to the
different $S^{N_j}M$ factors. Then it is clear from theorem \ref{t6.1}
that the moduli 
space of vortex solutions can be identified with the disjoint union of
all the $\Sigma_{N_0 , \ldots , N_n}(M)$ such that the non-negative
integers $N_j$ satisfy the condition $\alpha (P) = \sum_j \beta_j \,
N_j$ in ${\mathbb Z}^n$. By the argument after expression 
(\ref{7.2.7}) this condition is just $N_l - N_0 = \sum_k C_{lk}\, 
\alpha_k (P)$ for all $l= 1, \ldots , n$.
The topological energy of each vortex solution also reduces to
\[
T_{[\phi]}\ =\ b\cdot \alpha (P) \; +\; \frac{\pi}{n+1} \sum_{0\leq
  j\leq n} N_j  \ . 
\]
\end{rem}
\begin{rem}
It is manifest in theorem \ref{t6.1} that the vortex moduli space does
not change when the constant $c(P,M,a)$ or the moment map $\mu$ are
deformed by a translation, as long as this constant remains in the
interior or exterior of the image polytope $\mu (F) = \Delta$. 
When the constant $c(P,M,a)$ lies in the boundary of $\Delta$, then 
according to theorem \ref{t4.1} the solutions $(A, \phi)$ of the
vortex equations are constrained to satisfy $\phi (M) \subset E_{j_1}
\cap \cdots \cap E_{j_k}$, where $n-k$ is the dimension of the open
face of $\Delta$ that contains $c(P,M,a)$. Thus in some sense this
case corresponds to sigma-models with target $\CC {\mathbb P}^{n-k}$
and gauge group $T^n$. Since this gauge group is too big and has a
subtorus $T^k$ that acts trivially on the sections $\phi$, these cases
are somewhat degenerate.  
\end{rem}

\subsection{Proof of theorem \ref{t6.1} }

\subsubsection*{Equivalent theorem}

The first statement of theorem \ref{t6.1} follows from theorem
\ref{t4.1} and corollary \ref{c4.3}. As for the rest of  theorem
\ref{t6.1}, we will prove it  by stating and proving the equivalent 
theorem \ref{t6.3}.

Let ${\mathcal S}$ be the set of solutions of the vortex equations
(\ref{2.0}), and define
\[
{\mathcal B} \ =\ \left\{ (A,\phi) \in \AAA^{1,1}(P) \times \Gamma (E)
: \ \db^A \phi =0 \ \ {\rm and}\  \ \phi (M) \not\subseteq E_j \ \ {\rm for\ 
  all} \ 0\leq\ j\leq n \right\}\ .
\]  
The first thing to notice is that by theorem \ref{t4.1}, corollary \ref{c4.3}
and the subsequent remark, the natural inclusion of ${\mathcal S}$ in
${\mathcal B}$ actually induces a bijection of quotient spaces
\begin{gather}
{\mathcal S}/ ({\rm real\ gauge\ transf.}) \ \longleftrightarrow \
{\mathcal B} / ({\rm complex\ gauge\ transf.}) \ . 
\label{7.2.1}
\end{gather}
On the other hand the action of $T^n_\CC \simeq (\CC^\ast )^n $ on
$\CC {\mathbb P}^n$ is given by (\ref{6.1}), with $g_j \in \CC^\ast$,
so it is clear from the definitions of $F_j$ and $E_j$ that
\[
\phi^{-1} (E_j)\ =\ (g\cdot \phi )^{-1} (E_j) \qquad \subset \ \ M
\]
for any complex gauge transformation $g:M\rightarrow
T^n_\CC$. Moreover, it is a direct consequence of propositions \ref{p6.4}
and \ref{p6.5}, stated below, that for any irreducible hypersurface
$Z\subset M$ the intersection multiplicities satisfy
\[
{\rm mult}_{\phi (Z)} (E_j ,\phi (M)) \ =\ {\rm mult}_{(g\cdot
  \phi)(Z)} (E_j ,\; (g\cdot \phi) (M)) \ ,
\]
or in other words they are complex-gauge invariant. This fact
together with the bijection (\ref{7.2.1}) (which, recall, is induced by the
inclusion ${\mathcal S} \hookrightarrow {\mathcal B} $ ) show that
theorem \ref{t6.1} is equivalent to the following result.
\begin{thm}
Assume that $c(P,M,a)$ lies in the interior of $\Delta$, and for each
$j=0,\ldots , n$ pick an effective divisor $D_j = \sum_i a^i_j \cdot
Z_i$ on $M$ such that conditions (i) and (ii) of theorem \ref{t6.1} are
satisfied. Then there exists a pair $(A,\phi ) \in {\mathcal B}$,
unique up to complex gauge equivalence, such that 
\begin{gather}
{\rm mult}_{\phi (Z_i)} (E_j , \phi (M))\ =\ a^i_j \ .
\label{7.2.2}
\end{gather}
Furthermore all pairs in ${\mathcal B}$ can be obtained in this way.
\label{t6.3}
\end{thm} 

The method that we will use to prove this theorem is not intrinsic, in
the sense that it is based on the use of the usual local charts from
$\CC {\mathbb P}^n$ to $\CC^n$. In informal terms, we use the fact that
the domains of these charts are $T^n$-invariant and dense in $\CC {\mathbb
P}^n$ to transfer the problem of finding holomorphic sections of $E$
--- which has fibre $\CC {\mathbb P}^n$ --- to the problem of finding
meromorphic sections of vector bundles with fibre $\CC^n$.

\subsubsection*{Proof of the equivalent theorem}

As always, we start by introducing some notation. For each $j = 0,
\ldots , n$ define the action $\rho_j$ of $T^n$ on $\CC$ by
restricting the action $\rho$ of formula (\ref{6.1}) to the $j$-th
homogeneous coordinate of $\CC {\mathbb P}^n$. Thus for example $\rho_0$
is the trivial action, while for $j\ne 0$ the actions $\rho_j$ depend
on the matrix $C$. Define also the associated line bundles $L_j = P
\times_{\rho_j} \CC$.

Now consider the usual complex charts $\varphi_j : \UU_j \rightarrow
\CC^n$ of $\CC {\mathbb P}^n$, defined by 
\begin{align}
\UU_j & \ =\  \CC {\mathbb P}^n \setminus F_j \ =\ \left\{ [z_0 , \ldots ,
  z_n ] \in \CC {\mathbb P}^n : z_j \ne 0   \right\} \ ,   \nonumber   \\
\varphi_j & ([z_0 , \ldots, z_n]) \ =\ z_j^{-1} \ (z_0 , \ldots ,
  z_{j-1}, z_{j+1}, \ldots , z_n) \ .
\label{7.2.8}
\end{align}
From formula (\ref{6.1}) and the definition of $\rho_j$ it is clear
that $\UU_j$ is $T^n$-invariant and that, for any $g \in T^n$,
\[
\varphi_j \circ \rho_g \ =\ \Big(\; \ldots ,\;(\rho_{j})_{g^{-1}}
\cdot (\rho_{k})_g ,\; \ldots \; \Big)_{k\ne j} \circ \varphi_j \ . 
\]
Thus defining the vector bundle over $M$
\begin{gather}
V_j \ =\ (L_j)^{-1} \otimes \left( L_0 \oplus \cdots \oplus L_{j-1}
\oplus L_{j+1} \oplus \cdots \oplus L_n \right) \ ,
\label{7.2.3}
\end{gather}
and recalling that 
\[
E\setminus E_j \ =\ \big\{  [p,x] \in P\times_\rho \CC {\mathbb P}^n :\
x\in \UU_j  \big\}\ ,
\]
one has that the maps
\begin{gather}
\tvp_j : E\setminus E_j \rightarrow V_j \ , \qquad [p, x] \mapsto
    [p,\varphi_j (x)]
\label{7.2.9}
\end{gather}
are well defined. These maps clearly are fibre-preserving
diffeomorphisms. As a matter of notation, we will sometimes call
$L_{j,k}$ the $k$-th line bundle in the direct sum decomposition (\ref{7.2.3});
thus for example
\begin{gather}
L_{0,k}\ =\ (L_0)^{-1}\otimes L_k \ , \ \quad L_{n,k}\ =\ (L_n)^{-1} \otimes
L_{k-1} \ \quad {\rm and} \quad   V_j \ =\ \oplus_{1\leq k\leq n} \;
L_{j,k}\ .
\label{7.2.4}
\end{gather}
Since the actions $\rho_j$ preserve the canonical hermitian product on
$\CC$, the line bundles $L_{j,k}$ are all equipped with a natural
hermitian metric, denoted $h_{j,k}$. Another standard fact is that a
connection $A$ on $P$ induces connections on the associated line
bundles $L_j$ and $L_{j,k}$. These connections are
$h_{j,k}$-compatible. If the connection $A$ is integrable, i.e
$F_A^{0,2} =0$, then the induced connections on the $L_{j,k}$ are
integrable as well, i.e. their curvature form is in $\Omega^{1,1}(M)$.

The reason why we are interested in these integrable connections is
that, according to a well known result, an integrable, metric-compatible
connection $\nabla$ on a $C^{\infty}$ hermitian vector bundle
$(V,h)\rightarrow M$, induces a unique holomorphic structure
${\mathcal H}$ on $V$ such that $\nabla$ is the hermitian connection
of $(V,h,{\mathcal H})$ \cite{Ko}. The bundle $V$ equipped with this
holomorphic structure will be denoted by $V^{\nabla}$. We will often
apply this result to the line bundles $L_{j,k}$. When the integrable
connection on $L_{j,k}$ comes from a connection $A$ on $P$, we denote
by $L_{j,k}^A$ the line bundle together with the induced holomorphic
structure.  

Using all these conventions we define 
\[
{ \mathcal C } \ =\ \left\{ (\nabla_1 , \xi_1 , \ldots , \nabla_n ,
\xi_n )  : \  {\rm conditions\ (1)\ and\ (2)\ are\ satisfied}  \right\}
  \ ,
\]
where the conditions are 
\begin{itemize}
\item[(1)] $\nabla_k$ is an $h_{0,k}$-compatible connection on
  $L_{0,k}$ and  $\xi_k$ is a non-zero meromorphic section of
  $L^{\nabla_k}_{0,k}$; 
\item[(2)] the divisors on $M$ associated to the sections $\xi_k$
  satisfy $(\xi_1)_{-} = \cdots = (\xi_n)_{-} =: (\xi)_{-}$, and the
  intersection ${\rm supp}\, (\xi)_{-}\, \cap \: {\rm
  supp}\, (\xi_1)_{+} \,  
  \cap\, \cdots \,\cap \: {\rm supp}\, (\xi_n)_{+}$ is empty.
\end{itemize}
In the last condition we have decomposed a divisor $D = D_{+} - D_{-}$
into its positive and negative parts. The main tools to prove theorem
\ref{t6.3} are then the following two propositions.
\begin{prop}
There exists a bijection $\Upsilon:{\mathcal B} \rightarrow {\mathcal
  C}$ determined by the following conditions.
\begin{itemize}
\item[(i)] $\nabla_k$ is the connection on $L_{0,k}$ induced by the
  connection $A$ on $P$.
\item[(ii)] $\tvp_0 \circ \phi (q) \ =\ (\xi_1 (q) , \ldots , \xi_n
  (q))$ in $V_0$ for all $q\in M\setminus \phi^{-1}(E_0) $.
\end{itemize}
Furthermore, let $(A,\phi)$ and $(A' ,\phi')$ be two pairs in
${\mathcal B}$ and let $(\ldots , \nabla_k , \xi_k , \ldots )$ and
$(\ldots , \nabla_k' , \xi_k' , \ldots )$ be their images by
$\Upsilon$. Then $(A,\phi)$ and $(A' ,\phi')$ are complex gauge
equivalent if and only if for all $k$ the meromorphic sections $\xi_k$
and $\xi_k'$ have the same associated divisor in $M$.
\label{p6.4} 
\end{prop}

\begin{prop}
Given a pair $(A,\phi)$ in ${\mathcal B}$, let $(\ldots, \nabla_k ,
\xi_k, \ldots )$ be its image in ${\mathcal C}$ by
the bijection $\Upsilon$. Then for any irreducible analytic
hypersurface $Z \subset M$ and for any $j= 0,\ldots , n$, the
multiplicity of intersection of $E_j$ and $\phi (M)$ along $\phi (Z)$
is given by
\begin{gather}
{\rm mult}_{\phi (Z)} (E_j ,\phi (M)) \ = \ {\rm ord}_Z (\xi_j) -
\min_{0\leq k \leq n} \{ {\rm ord}_Z (\xi_k)\} \ ,
\label{7.2.5}
\end{gather}  
where we define ${\rm ord}_Z (\xi_0)\, =\, 0$.  
\label{p6.5}
\end{prop}

\begin{rem}
Formula (\ref{7.2.5}) implies that, for fixed $Z$, the multiplicities ${\rm
 mult}_{\phi (Z)} (E_j , \phi (M))$ are non-negative integers, with
 at least one of them being zero. On the other hand, given such a set
 of multiplicities, put
\begin{gather}
{\rm ord}_Z (\xi_j) \ =\ {\rm mult}_{\phi (Z)} (E_j , \phi (M)) - {\rm
  mult}_{\phi (Z)} (E_0 , \phi (M)) \ .
\label{7.2.6}
\end{gather}
It is then apparent that formulae (\ref{7.2.5}) and (\ref{7.2.6})
define inverse maps
between the set of sets of $n$ arbitrary integers ${\rm ord}_Z
(\xi_j)$, and the set of sets of $n+1$ non-negative, and not all
positive, integers.
\end{rem}
The proofs of these propositions, especially the first one, are rather
long and uninteresting, so will be exiled to appendix A. We are now ready
to prove theorem \ref{t6.3}. 

$\ $

\noindent
{\bf Proof of theorem \ref{t6.3}.}
Let $D_0 , \ldots , D_n$ be divisors on $M$ satisfying conditions
(i) and (ii) of theorem \ref{t6.1}. We will first show the
existence of a pair $(A, \phi) \in {\mathcal B}$ satisfying (\ref{7.2.2}).
By  definition of the actions $\rho_j$, the line bundle $L_0$ is
trivial and, for $k>0$,
\[
L_k \ =\ \bigotimes_{1 \le l \le n} (\hat{L}_l)^{C_{kl}} \ ,
\]
where the line bundles $\hat{L}_j$ were defined in section 2.4. In
particular, using (\ref{7.2.4}), this implies that
\begin{gather}
c_1 (L_{0,k})\ =\ \sum_l C_{kl} \  c_1 (\hat{L}_l) \ =\ \sum_l C_{kl} \;
\alpha_l (P)\ .
\label{7.2.7}
\end{gather}
Now denote by $e_1 , \ldots , e_n$ the standard basis of $\RR^n$ and
by $e_0$ the vector $-e_1- \cdots - e_n$. It is not difficult to check
directly that, for $a= 0, \ldots , n$, the vector $\beta_a = C^{-1}
e_a$ is a primitive vector in ${\mathbb Z}^n$ normal to one of the
$(n-1)$-dimensional faces of the polytope $\Delta$. Moreover this is
an outward pointing normal vector, and so the $\beta_a$'s coincide
with the $\beta_j$'s that appear in condition (ii) of theorem
\ref{t6.1}. Hence by this condition
\[
c_1 (L_{0,k}) \ =\ \sum_{a,l}\, C_{kl}\, (C^{-1} e_a)_l \ {\rm PD}(D_a)\
=\ {\rm PD}(D_k -D_0)\ .
\]
Thus the divisor $D_k - D_0$ defines a holomorphic structure on
$L_{0,k}$ together with a meromorphic section $\xi_k$ of this line
bundle (see \cite{Ko} or \cite{G-H}). Now denote by $\nabla_k$ the hermitian
connection on the hermitian bundle $(L_{0,k} ,\, h_{0,k} )$ equipped
with that holomorphic structure. Then by construction the multiplet
$(\nabla_1 , \xi_1 , \ldots , \nabla_n , \xi_n )$ satisfies condition
(1) of the definition of ${\mathcal C}$. Since $(\xi_k)_{+} = D_k $
and $(\xi_k)_{-} = D_0 $, it also satisfies condition (2), as follows
from the requirement (i) of theorem \ref{t6.1}. Thus
\[
( \ldots , \nabla_k , \xi_k , \ldots )_{1\le k \le n} \ \in \
  {\mathcal C}\ .
\]
According to propositions \ref{p6.4} and \ref{p6.5} this determines a
pair $(A, \phi) \in {\mathcal B}$ such that
\begin{gather*}
{\rm mult}_{\phi (Z_i)} \, (E_j , \phi (M)) \ =\ (a^i_j - a^i_0)\; -
\; \min_{0 \leq k \leq n} \{ a^i_k - a^i_0  \}\  
=\ a^i_j \; - \; \min_{0 \leq k \leq n} \{ a^i_k \} \ =\ a^i_j \ ,
\end{gather*}
where in the last equality we have used again the requirement (i) on
the divisors $D_k$. This settles the existence part of theorem \ref{t6.3}.

We will now prove the uniqueness statement. Keeping the same divisors
$D_j = \sum_i a^i_j \cdot Z_i$ as above, suppose that $(A' , \phi')
\in {\mathcal B}$ is another pair that satisfies (\ref{7.2.2}), and denote
by
\[
( \ldots , \nabla'_k , \xi'_k , \ldots ) \ \in \
  {\mathcal C}
\]
the image of this pair by the bijection $\Upsilon$. Since we are
assuming that the intersection divisors of $\phi (M)$ and $\phi' (M)$
with $E_j$ are the same, proposition \ref{p6.5} and formula (\ref{7.2.6})
in the subsequent remark imply that
\[
{\rm ord}_{Z_i} (\xi_j) \ =\ a^i_j - a^i_0 \ =\ {\rm ord}_{Z_i}
(\xi'_j)\ .
\]
Thus the meromorphic sections $\xi_j$ and $\xi'_j$ have the same
divisor in $M$, and from proposition \ref{p6.4} we conclude that $(A'
, \phi')$ is complex gauge equivalent to $(A ,\phi )$, as required.

To complete the proof of theorem \ref{t6.3} we just need to justify
the last assertion, i.e. that for every pair $(A,\phi) \in {\mathcal
  B}$ the divisors 
$D_j = \sum_i a^i_j \cdot Z_i$ defined by (\ref{7.2.2}) satisfy conditions
(i) and (ii) of theorem \ref{t6.1}. In the first place, 
the definition of ${\mathcal B}$ tells us that $\phi (M) \not\subseteq
E_j$, so $\phi^{-1} (E_j)$ is a union of irreducible hypersurfaces of
$M$. In particular the intersection multiplicities of (\ref{7.2.2}) are 
finite integers, and the divisors $D_j$ are well defined. Secondly, by
the definition of $D_j$ as the inverse image by $\phi$ of the 
intersection divisor of $\phi (M)$ and $E_j$, we have that $\phi ({\rm
  supp}\; D_j ) \subset E_j $. This implies that (i) is satisfied, because
$\cap_j E_j = \emptyset$. 
Finally, to recognize that the divisors $D_j$ associated to $\phi$
satisfy (ii) as 
well, consider the sections $\xi_k = (\tvp_0 \circ \phi )_k$ of 
proposition \ref{p6.4}. From proposition \ref{p6.5} and
formula (\ref{7.2.6}) it is clear that $D_k - D_0 $
is just the divisor of $\xi_k$. But $\xi_k$ is a meromorphic section
of $L^A_{0,k}$, and so by standard
results the Poincar\'e dual of the divisor of $\xi_k$ is $c_1
(L_{0,k})$. Using these facts, (\ref{7.2.7}), and the formula
$\beta_a = C^{-1} e_a$ established earlier, we then get that
\[
\sum_{0\le a \le n} C \,\beta_a \; {\rm PD}(D_a) \ =\ \sum_{1\le k \le
n} e_k \; {\rm PD}(D_k - D_0) \ =\ \sum_{1\le k \le n} e_k  \; c_1
(L_{0,k}) \ =\ C \, \alpha (P)
\]
in $H^2 (M; {\mathbb Z})^n$. Multiplying on the left by the matrix
$C^{-1}$ we obtain that (ii) is indeed satisfied. \hfill $\qed$ 
\vspace{0.3cm}

\subsection{Proof of theorem \ref{t6.2} }

The main task is to express the cohomology class $[\eta_E] \in H^2 (E;
\RR )$ in terms of the Poincar\'e dual of $[E_j] \in H_{2(m+n)-2} (E ;
{\mathbb Z})$. We start by noticing that, up to exact forms, any
closed 2-form on $M$ may be written as $s_1 \, \omega_M + \beta$, where
$s_1 \in \RR$, and $\beta \in \Omega^2 (M)$ is such that $\beta \wedge
\omega_M^{m-1} = 0$. This is a consequence of the Lefschetz
decomposition on K\"ahler manifolds \cite{G-H}. Furthermore it is apparent
from expression (\ref{7.4.4}) below that the class of $\eta_E (A)$, when
restricted to the fibres $E_x \simeq \CC {\mathbb P}^n$ of $E$, generates
the group $H^2 (E_x ; \RR) \simeq \RR$. It is then a consequence of
the Leray-Hirsch theorem \cite{B-T} that ${\rm PD}(E_j)$, and in fact any
element of $H^2 (E; \RR)$, is of the form
\[
{\rm PD}(E_j) \ =\ s_0 \: [\eta_E] + s_1 \: \pi_E^\ast [\omega_M] +
\pi_E^\ast [\beta] \ ,
\] 
where $s_0 , s_1 \in \RR$. Hence we have that
\[
s_0 \int_M \phi^\ast [\eta_E] \wedge \omega_M^{[m-1]} \ =\ \int_M
\phi^\ast {\rm PD}(E_j) \wedge \omega_M^{[m-1]} \ -\ s_1 \int_M m
\; \omega_M^{[m]} \ .
\]
Now, by well known properties of the Poincar\'e duality, the
restriction of ${\rm PD}(E_j)$ to $\phi (M)$ is just the Poincar\'e
dual in $\phi (M)$ of the intersection divisor of $E_j$ and $\phi (M)$
\footnote{I thank Dr. J.M. Woolf for explaining this to me.}. But
because of (\ref{7.2.2}) this divisor is just $\sum_i a^i_j \cdot \phi
(Z_i) = \phi_\ast (D_j)$, and 
since $\phi : M \rightarrow \phi (M)$ is a biholomorphism we obtain
that
\[
\phi^\ast \ {\rm PD}(E_j) \ =\ \phi^\ast \ {\rm PD}(\phi_\ast \, D_j
)\ =\ {\rm PD}(D_j) \qquad {\rm in}\ \ H^2 (M; {\mathbb Z}) \ .
\]
Thus
\begin{gather}
s_0 \int_M \phi^\ast [\eta_E] \wedge \omega_M^{[m-1]} \ =\ -s_1 \;
m \; ({\rm Vol}\ M) \ + \ \int_M {\rm PD}(D_j) \wedge
\omega_{M}^{[m-1]} \ . 
\label{7.4.1}
\end{gather}
The task now is to compute the constants $s_0$ and $s_1$.
Firstly  we remark that
\begin{gather}
\int_{E_j} \: \eta_E^{[n-1]} \wedge \pi_E^\ast \, \omega_M^{[m]} \ =\ 
\int_E \: \eta_E^{[n-1]} \wedge \pi_E^\ast \: \omega_M^{[m]} \wedge {\rm PD}
(E_j) \ =\ n
\; s_0 \int_E \eta_E^{[n]} \wedge \pi_E^\ast \, \omega_M^{[m]} \ .
\label{7.4.2}
\end{gather}
Also
\begin{align}
\int_{E_j} \ \eta_E^{[n]} \wedge \pi_E^\ast \, \omega_M^{[m-1]} \ &= \
\int_E \: \eta_E^{[n]} \wedge \pi_E^\ast \: \omega_M^{[m-1]} \wedge {\rm PD}
(E_j)\ =    \label{7.4.3}     \\
&=\ (n+1) s_0 \int_E \eta_E^{[n+1]} \wedge \pi_E^\ast \, \omega_M^{[m-1]}\ +
\ m \, s_1 \int_E \eta_E^{[n]} \wedge \pi_E^\ast \, \omega_M^{[m]} \ .
\nonumber 
\end{align}
The constants $s_0$ and $s_1$ are therefore determined by the value of
the integrals in (\ref{7.4.2}) and (\ref{7.4.3}).
To compute these integrals, recall from section 2.2 that $[\eta_E]$ is
the cohomology class of the closed 2-form $\eta_E (A)$ in $E$, and
that
\[
\chi^\ast \eta_E (A) \ =\ \omega_F - \dd (\mu , A) \qquad {\rm in} \ \
\Omega^2 (P \times F)\ .
\]
As in (\ref{2.11}), a local section $s : \UU \rightarrow \pi_P^{-1} (\UU )$ of
$P$ determines a trivialization $E|_\UU \simeq \UU \times F$, and it is
not difficult to check that with respect to this trivialization we
have
\begin{gather}
 \eta_E (A)\ |_\UU \ =\ \omega_F - \dd (\mu , s^\ast A) \qquad {\rm in} \ \
\Omega^2 (\UU \times F)\ .
\label{7.4.4}
\end{gather}
It follows that for any $k \in {\mathbb N}$
\[
\eta_E (A)^{[k]} \ =\ \omega_F^{[k]}\ -\ \omega_F^{[k-1]} \wedge \dd
(\mu , s^\ast A )\ +\ \cdots  \qquad {\rm in} \ \ \Omega^{2k} (\UU
\times F) \ ,
\]
and integrating along the fibre \cite{B-T} we get that
\[
(\pi_E)_\ast \; (\eta_E (A)^{[k]}) \ =\ \left\{
\begin{aligned}
&\ {\rm Vol}\, F \quad \qquad \ {\rm if}\ \ k=n  \\
& -\, F_A^l \int_F \mu_l \quad {\rm if}\ \ k=n+1  \ .
\end{aligned}
\right.
\]
Using the standard properties of the homomorphism $(\pi_E)_\ast$ (see
\cite{B-T}), we therefore have that
\begin{align}
& \int_E \eta_E^{[n]} \wedge \pi_E^\ast \, \omega_M^{[m]} \ =\ \int_M
(\pi_E)_\ast (\eta_E^{[n]}) \wedge \omega_M^{[m]} \ =\ ({\rm Vol}\
F)\, ({\rm Vol}\ M) \ , 
\label{7.4.5}   \\
& \int_E \eta_E^{[n+1]} \wedge \pi_E^\ast \omega_M^{[m-1]} \ = \ \int_M
(\pi_E)_\ast (\eta_E^{[n+1]}) \wedge \omega_M^{[m-1]} \ =\ - \left(
\int_F \mu_l \right) \; \int_M F_A^l \wedge \omega_M^{[m-1]} \ .
\label{7.4.6}
\end{align}
Now consider the inclusions $i_{F_j} : F_j \hookrightarrow F$ and $i_{E_j} :
E_j \hookrightarrow E$. Using the restriction $i_{F_j}^\ast \omega_F$ as a
K\"ahler form on $F_j$, and $\mu \circ i_{F_j}$ as a moment map for the
$T^n$-action on $F_j$, one can define the 2-form $\eta_{E_j} (A)$ on
$E_j = P \times_{T^n} F_j$, which is the analogue of $\eta_E (A)$ on
$E$. But
\[
\chi^\ast \; i_{E_j}^\ast\; \eta_E (A) \ =\ i_{F_j}^\ast \; \omega_F - \dd (\mu
\circ i_{F_j} ,\, A) \ =\ \chi^\ast \; \eta_{E_j}(A) \ ,
\]
and so $\eta_{E_j}(A)$ is just $i_{E_j}^\ast  \eta_E (A)$. Since
$\pi_{E_j} = \pi_E \circ i_{E_j}$ as well, by analogy with
(\ref{7.4.5}) and (\ref{7.4.6}) we have that
\begin{align*}
\int_{E_j} i_{E_j}^\ast \; ( \eta_E^{[n-1]} \wedge \pi_E^\ast \,
\omega_M^{[m]} ) \ &=\ \int_{E_j} 
\eta_{E_j}^{[n-1]} \wedge \pi_{E_j}^\ast \omega_M^{[m]} \ =\ ({\rm Vol}\
F_j)\, ({\rm Vol}\ M)  \\
\int_{E_j} i_{E_j}^\ast ( \eta_E^{[n]} \wedge \pi_E^\ast \,
\omega_M^{[m-1]}) \ &= \ \int_{E_j} 
\eta_{E_j}^{[n]} \wedge (\pi_{E_j})^\ast  \omega_M^{[m-1]} \ =\ - \left(
\int_{F_j} \mu_l \circ i_{F_j} \right)  \int_M F_A^l \wedge
\omega_M^{[m-1]}  . 
\end{align*}
From the value of these integrals it is straightforward to compute the
constants $s_0$ and $s_1$; it is enough to use (\ref{7.4.2}),
(\ref{7.4.3}) and the fact 
that $\CC {\mathbb P}^k$ with the Fubini-Study metric has volume $\pi^k /
k!$. Doing this and substituting the result into (\ref{7.4.1}), one obtains
that
\begin{eqnarray*}
\int_M \phi^\ast [\eta_E] \wedge \omega_M^{[m-1]}\ &=&\ \frac{n!}{\pi^n}
\left[ (n+1) \int_F \mu \ - \pi \int_{F_j} \mu \circ i_{F_j}  \right] \
\cdot\ \int_M - F_A \wedge \omega_M^{[m-1]}\ + \\
&  & \qquad \qquad \qquad \qquad \qquad \qquad \qquad \ 
 + \ \pi \int_M {\rm PD}(D_j) \wedge \omega_M^{[m-1]} \ .
\end{eqnarray*}
Now on the one hand, as mentioned in section 2.4, the cohomology
class of $-F_A$ is just $\alpha (P)$. On the other hand, since the
equality above is valid for all $j$, we may as well sum over $j$ and
divide by $n+1$. Doing this, applying lemma \ref{lb.1} below and using the
definitions (\ref{2.6}) and (\ref{2.8}), we obtain the formula of theorem
\ref{t6.2}.

\begin{lem}
Denoting by $b \in \RR^n$ the barycentre of the polytope $\Delta$, one
has that
\begin{gather}
\frac{n!}{\pi^n} \left[ (n+1) \int_F \mu \ - \frac{\pi}{n+1}
  \sum_{j=0}^{n} \int_{F_j} \mu \circ i_{F_j}  \right] \ =\ \frac{1}{{\rm
  Vol\ }\Delta} \int_{v\in \Delta} v \ =\ b \ .
\label{7.4.7}
\end{gather}
\label{lb.1}
\end{lem} 
\begin{prooff}
Instead of computing these integrals directly, using (\ref{6.mm}), we will
evaluate them using the Duistermaat-Heckman theorem (see for instance
\cite{MD-S, G, ACS}). Since $F = 
\CC {\mathbb P}^n$ is a toric manifold, the Duistermaat-Heckman
polynomial is piecewise a constant; with our conventions it is 1 in
the interior of $\Delta$ and 0 in the exterior. Therefore
\[
\int_F \mu \ =\ \int_{v\in \Delta} v \ ,
\]
where the integral on the right-hand side is taken with respect to the
Lebesgue measure in $\RR^n$. To evaluate the other integrals of lemma
\ref{lb.1}, take a $T^{n-1}$-action on $F_j \simeq {\mathbb CP}^{n-1}$ of the
same kind as (\ref{6.1}), and let $\mu_j : F_j \rightarrow \RR^{n-1}$ be a
moment map for it. Since $\mu \circ i_{F_j}$ is $T^{n-1}$-invariant, it is
clear that it can be written as
\[
\mu \circ i_{F_j} \ =\ S \circ \mu_j + {\rm const.} \ ,
\]
where $S: \RR^{n-1} \rightarrow \RR^n$ is some linear embedding. Thus
using again the Duistermaat-Heckman theorem we obtain that
\[
\int_{F_j} \mu \circ i_{F_j} \ =\ \int_{v \in \mu_j (F_j)} S(v) + {\rm
  const.} \ =\ \frac{{\rm Vol}\ \mu_j (F_j)}{{\rm Vol}\ \mu (F_j)}\
  \int_{v \in \mu (F_j)} v \ ,
\]
where the prefactor of the last term is just the inverse of the
determinant of $S$ as a linear map from $\RR^{n-1}$ to its
image. Finally a third application of the Duistermaat-Heckman theorem
shows that
\[
{\rm Vol}\ \Delta \ =\ {\rm Vol}\ F \ =\ \frac{\pi^n}{n!} \ =\
\frac{\pi}{n} \ {\rm Vol}\ F_j \ =\ \frac{\pi}{n}\ {\rm Vol}\ \mu_j
(F_j) \ .
\]
Hence the left-hand side of (\ref{7.4.7}) is equal to
\begin{gather}
\frac{n+1}{{\rm Vol}\ \Delta}\ \int_{v \in \Delta} v \ -\
\frac{n}{n+1}\ \sum_{j=0}^{n} \frac{1}{{\rm Vol}\ \mu (F_j)} \
\int_{v \in \mu (F_j)} v \ .
\label{7.4.8}
\end{gather}
Now notice that the first term in the expression above is $n+1$ times
the barycentre vector of $\Delta$, while the second term is
$n(n+1)^{-1}$ times the sum of the barycentre vectors of the faces
$\mu (F_j)$ of $\Delta$. But for a polytope $\Delta \subset
\RR^n$ with vertices $p_0 , \ldots , p_n$ the barycentre vector is
\[
b\ = \ \frac{1}{{\rm Vol}\ \Delta} \int_{v \in \Delta} v \ =\
\frac{1}{n+1}\; (p_0 + \cdots + p_n ) \qquad  \in \ \ \RR^n \ . 
\]
In particular, applying this expression to the faces $\mu (F_j)$ of
$\Delta$, we get that the last term of (\ref{7.4.8}) is just
\[
-\frac{n}{n+1}\ \sum_{j=0}^{n} \frac{1}{n} \ (p_0 + \cdots + \hat{p}_j
 + \cdots + p_n )\ =\ - \frac{n}{n+1} \ (p_0 +\cdots + p_n )\ =\ -n\,
 b \ .
\]
Substituting this into (\ref{7.4.8}) we get the required result.
\end{prooff}

\section{Constructing solutions on quotient targets}

\subsection{Induced solutions}

Consider the gauged $\sigma$-model determined by the data of section
2.1, and let $H$ be a closed normal subgroup of $G$. In informal
terms, the aim of this section is to compare the vortex equations
defined for target $F$ with $G$-action, and the vortex equations
defined for target $F/ H_\CC$ with $G/H$-action. The main result
obtained is theorem \ref{t3.6}. We point out that this operation of
quotient on the target is more delicate than, for example, the product
of targets. In particular, the second vortex equation does not
have any natural behaviour under these quotients, and so theorem
\ref{t3.6} only concerns the set $\hat{\mathcal S}$ of solutions 
of the first and third vortex equations. Nevertheless,
the results in this section are
still useful, because in section 4 we showed that, in the abelian
case, the quotients $\hat{\mathcal S} / \GG_\CC$ are often very
similar to the usual quotients ${\mathcal S}/ \GG$ of vortex
solutions. All this will eventually be  
used in section 7 to exhibit non-trivial solutions of the vortex
equations when the target $F$ is a K\"ahler toric manifold.

$\ $

We now formalize the problem. Denote by $\hh
\subset \g $ the Lie algebra of $H$, by $G' = G/H$ the quotient group,
and by $\g' \simeq \g /\hh$ the Lie algebra of $G'$. The invariant
inner product on $\g$ induces a splitting $\g = \hh \oplus \hh^\perp$
with associated projections $\pi_1 : \g \rightarrow \hh$ and $\pi_2 :
\g \rightarrow \hh^\perp$; it also induces natural identifications
$\g^\ast \simeq \g$, $\hh^\ast \simeq \hh$ and $(\g')^\ast \simeq \g' \simeq
\hh^\perp$. Using these identifications it is not difficult to check
that $\pi_1 \circ \mu : F \rightarrow \hh$ is a moment map for the
action of $H$ on $F$.

Suppose, moreover, that $H_\CC$ acts freely on $F$, and that there
exists an element $a \in \pi_1 \circ \mu (F)$ which is invariant by
the coadjoint action of $H$ on $\hh^\ast$. Then, by standard results
\cite{Ki},
the quotient $F /H_\CC$ is a K\"ahler manifold in a natural way. The
K\"ahler structure on $F' := F/H_\CC $ depends on the choice of $a$
and can be characterized as follows. The complex structure on $F'$ is
the only one such that the projection $\pi_F : F\rightarrow F'$ is
holomorphic; the symplectic form $\omega_{F'}$ on $F'$ is determined
by the condition
\[
i_{Z_a}^\ast \, \pi_{F}^\ast \, \omega_{F'}\ =\ i_{Z_a}^\ast \, \omega_F \ ,
\]   
where $Z_a$ is the inverse image $(\pi_1 \circ \mu )^{-1} (a)$, and
$i_{Z_a}$ is the inclusion $Z_a \hookrightarrow F$. Note that it can
be shown that 
$Z_a$ is a $H$-invariant submanifold of $F$, and that $F' = Z_a / H $
\cite{Ki}. 
\vspace{.3cm}
\begin{rem}
When $F$ is compact it is never possible to find hamiltonian
$H$-actions such that $H_\CC$ acts freely. On the other hand, denoting
by $\mu_H$ the moment map of the $H$-action, there is a canonical
choice of an ${\rm Ad}_H^\ast$-invariant element $a \in \mu_H (F)$,
which is $a= \int_{x\in F} \mu_H (x)$. It can then be shown that if
$H$ acts freely on $\mu_H^{-1} (a)$, then $H_\CC \cdot \mu_H^{-1}(a)$
is an open subset of $F$ where the action of $H_\CC$ is free
\cite{Ki}. 
\end{rem}
\vspace{.3cm}
The group $G'$ acts naturally on $F'$ by the rule 
\begin{gather}
\pi_G (g) \cdot \pi_F (p) \ =\ \pi_F (g\cdot p) \quad \qquad \forall\ g\in
G,\ p\in F \ ,
\label{3.3.1}
\end{gather}
where $\pi_G : G \rightarrow G'$ is the quotient map. It is not
difficult to check that this is still a holomorphic hamiltonian
action. In fact, a moment map $\mu' : F' \rightarrow \g'^\ast \simeq
\hh^\perp$ for this action is determined by the formula
\begin{gather}
i_{Z_a}^\ast \, \pi_{F}^\ast \, \mu' \ =\ i_{Z_a}^\ast (\pi_2 \circ \mu)\ .
\label{3.3.2}
\end{gather}
Besides acting on $F$, the subgroup $H$ also acts freely on the
principal bundle $P$. Let $P'$ be the quotient space $P/H$ and let
$\zeta: P \rightarrow P'$ be the quotient map. The
group $G'$ acts naturally and freely on $P'$, and if we define
the projection $\pi_{P'} : P' \rightarrow M$ by  
\[
\pi_{P'} \ \circ \ \zeta \ = \ \pi_P \ ,
\]
it is apparent that $P'$ is the total space of a $G'$-bundle over $M$.
If $A \in \Omega^1 (P ; \g )$ is a connection on $P$, it is clear that
$(\dd \pi_G) \circ A$ descends to a form $A' \in \Omega^1 (P' ;
\g')$. This is a connection form on the bundle $P'$ \cite[p. 79]{K-N}.
Using that $P \times_{{\rm Ad}_G} \g' \simeq P' \times_{{\rm Ad}_{G'}}
\g'$, the
curvature form of $A'$ is $\dd \pi_G \circ F_A \, \in \Omega^2 (M; P
\times_{{\rm Ad}_G} \g')$, where $F_A \, \in \Omega^2 (M; P
\times_{{\rm Ad}_G} \g)$ is the curvature form of $A$. Identifying
$\g' \simeq \hh^\perp$ this curvature form is just $\pi_2 \circ
F_A$. In particular
\begin{gather}
F_A^{0,2}\ =\ 0 \quad \Rightarrow \quad F_{A'}^{0,2}\ =\ \pi_2 \circ
F_A^{0,2}\ =\ 0\ .
\end{gather}
Consider now the associated bundle $E' = P \times_{{\rm Ad}_G} F' = P'
\times_{{\rm Ad}_{G'}} F'$. There is a natural bundle map $E
\rightarrow E'$ determined by the formula 
\begin{gather}
[p,q] \ \mapsto \ [p, \pi_F (q)] \quad \qquad \forall\ p\in P,\ q\in F \ .
\label{3.3.3}
\end{gather}
As always, this induces a map on the space of sections 
\[
\Gamma (E)\  \longrightarrow\ \Gamma (E'), \qquad \phi \mapsto \phi' \ .
\]
Using the definition (\ref{2.10}) and the holomorphy of $\pi_F$ it is
then not difficult to check that 
\begin{gather}
\db^A \phi \ =\ 0 \quad \Rightarrow \quad \db^{A'} \phi' \ =\ 0\ .
\label{3.3.4}
\end{gather}
Hence, in terms of the spaces of solutions 
\begin{align*}
{\mathcal S}(P,E)\ &=\ \left\{  (A, \phi) \in \AAA (P) \times \Gamma (E): \
{\rm equations\ (\ref{2.0})\ are\ satisfied\,} \right\} \quad  \qquad {\rm
  and} \\
\hat{{\mathcal S}}(P,E)\ &=\ \left\{  (A, \phi) \in \AAA (P) \times
\Gamma (E): \ {\rm equations\ (\ref{2.0a})\ and\ (\ref{2.0c})\ are\
  satisfied \,} \right\}\  , 
\end{align*}
we have that the correspondence $(A, \phi) \mapsto (A' , \phi')$
defines a map 
\begin{gather}
\Upsilon\ :\ \hat{{\mathcal S}}(P,E)\ \longrightarrow\ \hat{{\mathcal S}}(P'
,E') \ .   
\label{3.3.5}
\end{gather}
We will now see how $\Upsilon$ behaves when we quotient by complex gauge 
transformations. 

$\ $

To start with, recall that the quotient map $\pi_G$ can be extended to
a homomorphism $(\pi_G)_\CC : G_\CC \rightarrow G'_\CC$, and that this
homomorphism induces an identification $G'_\CC \simeq G_\CC / H_\CC$. 
The homomorphism $(\pi_G)_\CC$ then defines a natural bundle map
\begin{gather}
P \times_{{\rm Ad}_G} G_\CC \longrightarrow P \times_{{\rm Ad}_G} G'_\CC \
, \qquad [p,g] \mapsto [\, p, (\pi_G)_\CC (g)\,]\ .
\label{3.3.6}
\end{gather}
As always, composition with this bundle map defines a map of sections
\[
\pi_{\GG_\CC}: \GG_\CC \ \longrightarrow \GG'_\CC \ , \qquad g \mapsto g'\ .
\]
This map clearly is a homomorphism of gauge groups. One can also check
that it has the following naturality property.
\begin{lem}
Let $(A, \phi)$ be any pair in $\AAA (P) \times \Gamma (E)$ and $g$
any gauge transformation in $\GG_\CC$. Then $g(\phi)' = g'(\phi')$ in
$\Gamma (E')$ and $g(A)' = g' (A')$ in $\AAA (P')$.
\label{l3.5}
\end{lem}

A direct consequence of this lemma is that the map $\Upsilon$
descends to a map of quotient spaces $\hat{{\mathcal S}}(P,E) / \GG_\CC
\rightarrow \hat{{\mathcal S}} (P',E') / \GG'_\CC$. Another important
property of $\Upsilon$ is the following.
\begin{lem}
Let $(A_1 , \phi_1)$ and $(A_2, \phi_2)$ be two pairs in
 $\hat{{\mathcal S}}(P, E)$. Then $(A'_1 , \phi'_1) = (A'_2, \phi'_2)$ if
 and only if  
 there exists a gauge transformation $g \in {\mathcal H}_\CC $ such
 that $\phi_2 = g (\phi_1)$ and $A_2 = g(A_1)$. When it exists, this
 transformation is unique.
\label{p3.8}
\end{lem}
This is proved using the assumptions that $H_\CC$ acts
freely on $F$, the rules (\ref{3.1.1}) and (\ref{3.1.2}) for gauge
transformations, and the fact that the $(A_j , \phi_j )$ satisfy the
vortex equations (\ref{2.0a}) and (\ref{2.0c}). The details are exiled
to section 6.2. Combining the two lemmas above one directly obtains the main
result of this section, which is the following. 
\begin{thm}
The induced map $\Upsilon : \hat{{\mathcal S}}(P,E)/ \GG_\CC
\longrightarrow \hat{{\mathcal S}}(P',E') / \pi_{\GG_\CC} (\GG_\CC )$
is injective. 
\label{t3.6}
\end{thm}
\begin{rem}
In many cases of interest the homomorphism $\pi_{\GG_\CC}$ is
surjective, and so we actually get an injection $\Upsilon : \hat{{\mathcal
 S}}(P,E)/ \GG_\CC \rightarrow \hat{{\mathcal S}}(P',E') /
\GG'_\CC$. This happens, for 
example, when the group $G$ can be factorized as $G= H \times W$,
where $W$ is some other subgroup of $G$. In this case $G_\CC \simeq
H_\CC \times W_\CC$ and $G'_\CC \simeq W_\CC$, and so it is clear that
any section of $P \times_{{\rm Ad}_G} G'_\CC$ can be lifted to a
section of $P \times_{{\rm Ad}_G} G_\CC$.
\end{rem}

\subsection{Proof of lemma \ref{p3.8}}

We first establish a preparatory lemma, and then prove Lemma \ref{p3.8}.

\vspace{.2cm}
The proof of theorem \ref{t3.6} is essentially contained in proposition
\ref{p3.8} below. We state and prove this proposition after
establishing a preparatory 
lemma. Finally a few lines are spent completing the proof of theorem
\ref{t3.6}.

\begin{lem}
Let $\phi_1$ and $\phi_2$ be two sections of $E$. Then $\phi'_1 =
\phi'_2$ in $\Gamma (E')$ if and only if there exists a gauge
transformation $g \in {\mathcal H}_\CC $ such that $\phi_2 = g
(\phi_1)$. When it exists, this transformation is unique.
\label{l3.7}
\end{lem}
\begin{prooff}
Since $H_\CC$ acts freely on $F$, the gauge group ${\mathcal H}_\CC$
acts freely on 
$\Gamma (E)$, and so the uniqueness of $g$ is clear. Furthermore, the
sufficiency part is apparent from the formula $g(\phi)' = g' (\phi')$
of lemma \ref{l3.5},
so we just need to prove the necessity.

Suppose then that $\phi'_1 = \phi'_2$, and take local trivializations
of $E = P \times_G F$ and $E' = P \times_G F'$ induced by the same
local trivialization of $P$. With respect to these trivializations, if
the $\phi_i$ are represented by local maps $\hat{\phi}_i : \UU
\rightarrow F$, the $\phi'_i$ are represented by $\pi_F \circ
\hat{\phi}_i : \UU\rightarrow F' $. But the equality of the
$\phi'_i$'s implies that $\pi_F \circ \hat{\phi}_1 = \pi_F \circ
\hat{\phi}_2 $, and since $\pi_F : F \rightarrow F'$ is a principal
$H_\CC$-bundle we conclude that there exists a unique smooth map
$\hat{g} : \UU \rightarrow H_\CC $ such that 
\begin{gather} 
\hat{\phi}_2 (x)\ =\ \hat{g}(x)\, \cdot\, \hat{\phi}_1 (x) \qquad {\rm
  for\ all}\ x\in \UU \ .
\label{3.3.9}
\end{gather}
Had we chosen initially a different trivialization of $P$, related to
the first one by a transition function $\alpha : \UU \rightarrow G$,
we would get the maps $\alpha (\cdot)^{-1} \cdot \hat{\phi}_i (\cdot ) :
\UU \rightarrow F$ as representatives of the $\phi_i$'s. In particular
\[
\alpha (x)^{-1}\, \cdot\, \hat{\phi}_2 (x) \ =\ \left( \alpha (x)^{-1}
\hat{g}(x)\, \alpha (x) \right) \, \cdot \, \alpha (x)^{-1} \hat{\phi}_1
(x) \qquad {\rm for\ all}\ x\in \UU \ ,
\]
which shows that the local maps $\hat{g} (\cdot)$ transform as
sections of the bundle $P \times_{{\rm Ad}_G} H_\CC$. By their local
uniqueness, these maps can be ``glued'' together to define a global
section of the latter bundle, i.e. an element $g$ of ${\mathcal
  H}_\CC$. It then follows from (\ref{3.3.9}) that
$\phi_2 = g  (\phi_1)$.  
\end{prooff}

\noindent
{\bf Proof of Lemma \ref{p3.8} .}
By the previous lemma, the condition $\phi'_1 = \phi'_2$ is equivalent
to the existence of a unique $g \in {\mathcal H}_\CC$ such that
$\phi_2 = g (\phi_1)$. We now have to show that, for this
transformation $g$, $A'_2 = A'_1$ if and only if $A_2 = g (A_1)$.

A first observation is that, since $H$ is normal in $G$, for any $h
\in H$ and $v \in \g$ the vector ${\rm Ad}_h (v) - v$ is in
$\hh$. It is then apparent from the gauge transformation rule 
\[
g(A) \ =\ {\rm Ad}_g \circ A \ -\ \pi_P^\ast (g^{-1} \db g + \bar{g}^{-1}
\partial \bar{g}) \ , \qquad g\in {\mathcal H}_\CC \ ,
\]
that, for any connection $A$ on $P$, the difference $g(A) - A$ is in
$\Omega^1 (P ; \hh)$. Therefore, if $A_2 = g (A_1)$, using the
definition of $A'_j$ in the paragraph after (\ref{3.3.2}), we have that 
\[
\zeta^\ast (A'_2 - A'_1) \ =\ (\dd\pi_G) (A_2 -A_1) \ =\ 0 \ ,
\]
which implies that $A'_2 = A'_1$.
Conversely, suppose that $A'_2 = A'_1$. Then by the previous formula
$A_2 -A_1 \in \Omega^1 (P ; \hh)$, and therefore
\begin{gather}
g(A_1) - A_2 \ =\ g(A_1) - A_1 + A_1 - A_2 \quad \in \ \Omega^1 (P; \hh).
\label{3.3.10}
\end{gather}
On the other hand the assumptions of the proposition are that $(A_i ,
\phi_i) \in \hat{{\mathcal S}}(P, E)$, and since equation (\ref{2.0a})
is invariant by complex gauge transformations,
\begin{gather}
\db^{A_2} \phi_2 \ =\ \db^{A_1} \phi_1 \ =\ \db^{g(A_1)} g(\phi_1)\ =\
\db^{g(A_1)} \phi_2 \ =\ 0 \ .
\label{3.3.11}
\end{gather}
Now let $s: \UU \rightarrow P$ be a local trivialization, and let
$\hat{\phi}_i : \UU \rightarrow F$ be the representative of $\phi_i$
with respect to the associated trivialization of $E = P\times_G
F$. From (\ref{3.3.11}) and (\ref{2.12}) we get that
\[
\left[ \; (s^\ast A_2)_x  (\cdot) \; \right]^\flat \; |_{\hat{\phi}_2 (x)} \
=\ \left[ \; (s^\ast g(A_1))_x  (\cdot) \; \right]^\flat \; |_{\hat{\phi}_2
  (x)} \qquad {\rm for\ all}\ x\in \UU \ .
\]
Since the correspondence $\flat : \g \rightarrow T_q F, \ v \mapsto
v^\flat |_q$ is linear, this is the same as
\begin{gather}
\left[\ s^\ast (A_2 - g(A_1))_x \; (\cdot) \ \right]^\flat \; |_{\hat{\phi}_2
  (x)} \ =\ 0 \qquad  {\rm for\ all}\ x\in \UU\ . 
\label{3.3.12}
\end{gather}
But by assumption the action of $H$ on $F$ is free, and so the
restriction of $\flat$ to $\hh \subset \g$ is injective. From
(\ref{3.3.10}) and  (\ref{3.3.12}) it then follows that $s^\ast (A_2 -
g(A_1)) = 0$, and since the 
trivialization $s$ was arbitrary we conclude that $A_2 = g
(A_1)$.     \qed

\section {Solutions for target a compact toric  manifold }

\subsection{The canonical K\"ahler toric manifolds}

A compact K\"ahler toric manifold $F$ is by definition a compact
K\"ahler manifold equipped with an effective hamiltonian action of
$T^n$ --- where $n$ is the complex 
dimension of $F$ --- which operates by holomorphic transformations. If
$\mu : F \rightarrow \RR^n$ is  a moment map for this action, it is
well known that the image $\mu (F)$ is a special kind of polytope in
$\RR^n$, usually called a Delzant polytope, and that this polytope
determines $F$ up to $T^n$-equivariant symplectomorphisms \cite{Del}.
\begin{defn}
A Delzant polytope $\Delta$ in $\RR^n$ is a convex polytope such that:
\begin{itemize}
\item[$\bullet$] there are $n$ edges meeting at each vertex;
\item[$\bullet$] the edges meeting at the vertex $p$ are rational, in
  the sense that they are of the form $p + t v_i, \ t\in \RR ,$ with
  $v_i \in {\mathbb Z}^n$;
\item[$\bullet$] these $v_1 , \ldots , v_n$ can be chosen to be a
  basis of ${\mathbb Z}^n$.
\end{itemize}
\end{defn}
The symplectomorphism mentioned above, however, does not necessarily
preserve the complex structure on $F$, and so the polytope $\mu (F)$
does not determine $F$ as a K\"ahler manifold. In other words, this
means that several inequivalent K\"ahler toric manifolds may give rise
to the same image polytope $\mu (F)$. Although lacking injectivity,
the correspondence between K\"ahler toric manifolds and Delzant
polytopes is certainly surjective. This is because, given any Delzant
polytope $\Delta$ in $\RR^n$, there is a natural way to construct a
K\"ahler toric manifold $F_\Delta$ such that $\mu (F_\Delta) = \Delta
$. We will now briefly recall this construction; for more details see
for example \cite{G} or \cite{ACS}.

$\ $

Let $\Delta$ be a Delzant polytope in $\RR^n$ with $(n-1)$-dimensional
faces, or facets, $B_1 , \ldots , B_d$, where $d>n$. Then one can
uniquely choose 
vectors $u_1 , \ldots , u_d \, \in {\mathbb Z}^n$ such that $u_i$ is a
primitive, outward pointing, normal vector to $B_i$. The polytope
$\Delta$ is then the intersection  of half-spaces
\[
\left\{ x\in \RR^n : \ u_i \cdot x \ \le \ \lambda_i \ ,\ i=1,\ldots
,d  \right\} \ ,
\]
for some $\lambda_i \in \RR$. Denoting by $e_1 , \ldots ,e_d$ the
standard basis of $\RR^d$, define the linear map
\[
\beta : \RR^d \longrightarrow \RR^n \ , \qquad e_j \mapsto u_j \ ,
\]
and its $i$-linear extension $\beta_\CC : \CC^d \rightarrow \CC^n $.
It is not difficult to show that $\beta ({\mathbb Z}^d) = {\mathbb
  Z}^n$, and so these maps descend to homomorphisms of tori
\[
\begin{CD}
\RR^d  @>{\beta}>> \RR^n  @.\qquad \qquad @.  \CC^d  @>{\beta_\CC}>> \CC^n  \\
@VVV               @VVV         @.                  @VVV          @VVV   \\
T^d  @>{\tilde{\beta}}>> T^n  @. \qquad \qquad @.  T^d_\CC @>{\tilde{\beta}_\CC}>> T^n_\CC
 \ .
\end{CD}
\]
In both these diagrams the vertical arrows represent the exponential
map (\ref{5.1.1}). The subspace ${\mathfrak n}= \ker{\beta}$ of $\RR^d$
exponentiates to the subgroups
\begin{align*}
N\ &=\ \ker{\tilde{\beta}} \ =\ \exp{({\mathfrak n})}  \qquad \qquad
\qquad \quad \  \subset  \ T^d  \\
N_\CC \ &=\ \ker{\tilde{\beta}_\CC}\ =\ \exp{({\mathfrak n})} \times
\exp{(i{\mathfrak n})}  \qquad \ \subset T^d_\CC \ ,
\end{align*}
and one has the short exact sequence
\begin{gather}
0 \longrightarrow N_\CC \longrightarrow T^d_\CC  \longrightarrow
T^n_\CC \longrightarrow 0 \ .
\label{6.1.1}
\end{gather}
Now consider the natural action of $T^d$ on the K\"ahler manifold
$\CC^d$ given by 
\begin{gather}
(g_1 , \ldots , g_d) \cdot (z_1 , \ldots , z_d ) \ =\ (g_1 z_1 ,
  \ldots , g_d z_d )\ .
\label{6.1.2}
\end{gather}
This action operates by holomorphic transformations and has moment map
\[
\mu : \CC^d \longrightarrow \RR^d \ , \qquad (z_1 ,\ldots , z_d)
\ \mapsto\ -\pi (|z_1|^2 , \ldots , |z_d|^2 ) + (\lambda_1 , \ldots ,
\lambda_d) \ . 
\]
The restriction of this action to the subgroup $N$ has moment map
$\pi_1 \circ \mu : \CC^d \rightarrow {\mathfrak n}$, where $\pi_1$ is
the orthogonal projection from $\RR^d$ to ${\mathfrak n}$. Notice also that
the action (\ref{6.1.2}) has a natural extension to the complexified group
$T^d_\CC$; this is given by the same formula, but with the $g_j$'s
belonging to $\CC^\ast$. Now define the subset
\begin{gather*} 
\CC^d_\Delta \ =\ \left\{ (z_1 , \ldots , z_d )\in \CC^d : z_{j_1}=\cdots = z_{j_k} = 0\  
   {\rm is\ allowed\ only \ if}\ \bigcap_{1 \le l \le k} B_{j_l} \ne \emptyset
   \right\} \ .
\label{6.1.3}
\end{gather*}
It is shown in appendix 1 of \cite{G} that $\CC^d_\Delta$ is an open
dense subset 
of $\CC^d$ where $N_\CC$ acts freely. Furthermore the inverse image $Z
= (\pi_1 \circ \mu )^{-1} (0)$ is contained in $\CC^d_\Delta$, and in
fact $\CC^d_\Delta = N_\CC \cdot Z$. Hence, by the quotient
construction described in section 6, the quotient manifold $F_\Delta
= \CC^d_\Delta / N_\CC$ has a unique structure of K\"ahler manifold
such that the projection $\pi_F : \CC^d_\Delta \rightarrow F_\Delta$
is holomorphic and
\[
i_Z^\ast \; \pi_F^\ast \; \omega_{F_\Delta} \ =\ i_Z^\ast \; \omega_{\CC^d}\ .
\]
Just as in section 6, the quotient group $T^d /N $ acts
naturally on $F_\Delta$ by holomorphic transformations. Identifying
$T^d /N \simeq T^n $ through $\tilde{\beta}$, this action has a moment
map $\mu' : F_\Delta \rightarrow \RR^n$ determined by 
\[
i_Z^\ast \; \pi_F^\ast  \; \mu' \ =\ i_Z^\ast \; (\beta \circ \mu )\ .
\]
It can be shown that $\mu' (F_\Delta) = \Delta$. The K\"ahler manifold
$F_\Delta$ equipped with this $T^n$-action is the canonical K\"ahler
toric manifold we were looking for.

\begin{exa}
When $\Delta$ is the Delzant polytope
\[
\Delta \ =\ \left\{ x\in \RR^n : x_j \le 0 \ \ {\rm and}\ \ \Sigma_j \ x_j \ge
- \pi \right\} \ ,
\]
one gets that $\CC^d_\Delta = \CC^{n+1} \setminus \{ 0 \}$ and that $N
\simeq T^1$ is the diagonal subgroup of $T^d = T^{n+1}$. It is then
clear that $F_\Delta = \CC^d_\Delta / N_\CC = \CC {\mathbb P}^n$, and one
can check that the induced K\"ahler metric on $\CC {\mathbb P}^n$ is the
Fubini-Study one.
\end{exa}

\vspace{.2cm}

Besides the result $\mu' (F_\Delta) = \Delta$ described above, in the
next section we will also use that for any facet $B_j$ of $\Delta$
\begin{gather}
(\mu' \circ \pi_F)^{-1} (B_j) \ =\ \{ z\in \CC^d_\Delta : z_j =0   \}\ .
\label{6.1.4}
\end{gather}
This fact also follows from the results in the appendix 1 of \cite{G}.

\subsection{A family of non-trivial solutions}

Let $F$ be a compact K\"ahler toric manifold, let $\mu : F \rightarrow
\RR^n $ be a moment map for the associated torus action, and call
$\Delta$ the image $\mu (F)$, which is a Delzant polytope in $\RR^n$. 
Denote by $B_1 , \ldots , B_d$ the $(n-1)$-dimensional faces of
$\Delta$, and by $\beta_j \in {\mathbb Z}^n$ the unique primitive,
outward pointing, normal vector to $B_j$. 
Finally identify $\ttt^n \simeq \RR^n$ through
(\ref{5.1.1}), and take the euclidean inner product on $\RR^n$ to identify
$\ttt^n \simeq (\ttt^n)^\ast$.

$\ $

Now take any principal $T^n$-bundle $P'$ over the
K\"ahler manifold $M$, and denote by $\alpha (P') \in H^2 (M;{\mathbb
  Z})^n$ the vector of cohomology classes described in section
2.4. Using this principal bundle one can define the associated bundle
$E' = P' \times_{T^n} F $, which has base $M$ and typical fibre
$F$. From lemma \ref{l4.4} we know that the subsets
\[
F_j \ :=\  (\mu')^{-1} (B_j)
\]
are $T^n$-invariant complex submanifolds of $F$. Furthermore,
as described in the proof of lemma \ref{l4.5}, the associated bundles
\[
E'_j \ :=\ P' \times_{T^n} F_j
\]
are complex submanifolds of $(E' , J(A))$, where $J(A)$ is the complex
structure on $E'$ induced by an integrable connection $A$ on $P'$ (see section
2.2). The aim of this section is to prove the following result.

\begin{thm}
In the setting described above the vortex equations (\ref{2.0}) have
solutions only if the constant $c (P' , M, a)$ is in
$\Delta$. When this constant lies in the interior of $\Delta$, a set of
solutions can be described as follows. For each $j =1, \ldots , d$
pick an effective divisor $D_j = \sum_i a^i_j \cdot Z_i$ on $M$ such
that
\begin{itemize}
\item[$(i)$] if $\; \bigcap_{1 \le l \le k} B_{j_l} = \emptyset$ for some
  indices $j_1 , \ldots , j_k$, then the intersection of hypersurfaces
  ${\rm supp}\, D_{j_1} \cap \cdots \cap {\rm supp}\, D_{j_k}$ is empty
  as well;
\item[$(ii)$] the Poincar\'e duals of the fundamental homology cycles
  carried by the divisors $D_j$ satisfy $\alpha (P') = \sum_j \beta_j
 \; {\rm  PD}(D_j)$ in $H^2 (M; {\mathbb Z})^n$.
\end{itemize}
Then there is a solution $(A, \phi) \in {\mathcal S}(P', E')$ of the
vortex equations such that the intersection multiplicities of the
complex submanifolds $\phi (M)$ and $E'_j$ satisfy
\begin{gather}
{\rm mult}_{\phi (Z_i)} (E'_j , \phi (M)) \ =\ a^i_j \ .
\label{6.2.0}
\end{gather}
Different choices of divisors provide gauge inequivalent solutions.
\label{t5.3}
\end{thm}

Comparing with theorem \ref{t6.1} one recognizes that, when $F = \CC
{\mathbb P}^n$, the set of solutions obtained in theorem \ref{t5.3}
actually coincides with the full set of solutions, up to gauge
transformations. This motivates the following question.
\begin{ques}
Let F be any compact K\"ahler toric manifold, and suppose that the
constant $c(P' , M, a)$ lies in the interior of $\Delta$. Do the 
solutions described in theorem \ref{t5.3} represent the full space of
vortex solutions, up to gauge equivalence ?
\end{ques}

\noindent
{\bf Proof of theorem \ref{t5.3}.}
We first prove the theorem in the case where $F$ is the canonical
manifold $F_\Delta$. At the end we will deal with the case of any $F$
such that $\mu (F) = \Delta$.

Given the divisors $D_j$, by proposition \ref{p5.1} there is a principal
$T^d$-bundle $P \rightarrow M$ such that ${\rm PD}(D_j) = \alpha_j (P)
= c_1 (L_j) $ for all $j=1, \ldots , d$. Let
\begin{gather}
E \ =\ P\times_{T^d} \CC^d \ =\ L_1 \oplus \cdots \oplus L_d
\label{6.2.1}
\end{gather}
be the associated bundle. Using the notation of section 7.1, define
also the bundle $\dot{E} = P \times_{T^d} \CC^d_\Delta$, which is an
open dense subset of $E$.

Now consider the spaces of solutions ${\mathcal S}(P, E)$ and
$\hat{{\mathcal S}}(P, E)$ defined before (\ref{3.3.5}).
As in section 3 (with $C=Id$), since ${\rm PD}(D_j) = c_1 (L_j)$,
there exists a pair 
$(A, \phi) \in \hat{{\mathcal S}}(P, E)$ such that $D_j$ is the
divisor of the zero set of $\phi_j$ --- the $j$-th component of $\phi$
under the decomposition (\ref{6.2.1}) --- regarded as a holomorphic section of
$L^A_j$. Notice that condition (i) on the divisors implies that if
$\bigcap_{1 \le l \le k} B_{j_l} = \emptyset$ the intersection
$\bigcap_{1 \le l \le k} \phi^{-1}_{j_l} (0)$ is empty; thus, having in
mind the definition of $\CC^d_\Delta$, we conclude that $\phi (M)
\subset \dot{E}$, and the pair $(A, \phi)$ may be regarded as belonging
to $\hat{{\mathcal S}}(P, \dot{E})$.

At this point we want to apply the results of section 6  in order to
obtain solutions in $\hat{{\mathcal S}}(P', E')$. Going back to this
section, put $G = T^d, \ H = N$ and $F = \CC^d_\Delta$. The
homomorphism $\tilde{\beta} : T^d \rightarrow T^n$, which has kernel
$N$, provides identifications $T^d /N \simeq T^n$ and $P/N \simeq P
\times_{\tilde{\beta}} T^n$. But by lemma \ref{p5.2} and condition
(ii) we have that
\[
\alpha_a (P \times_{\tilde{\beta}} T^n) \ =\ \sum_l \beta_{al} \; 
\alpha_l (P) \ =\ \alpha_a (P') \qquad {\rm for\ all}\ a=1, \ldots , n.
\]
Thus the bundles $P/N$ and $P'$ are isomorphic, and so the $P'$ of
section 6 coincides with the $P'$ of this section. On the other hand,
since $F_\Delta$ is the K\"ahler quotient $\CC^d_\Delta / N_\CC$ with
the $T^n$-action provided by the identification of $T^d /N$ with $T^n$
through $\tilde{\beta}$,
the $F'$ and $E'$ of section 6 are just the $F_\Delta$ and $E'$ of
this section. Applying the results of section 6  we therefore have
that the map $\Upsilon$ takes $(A, \phi)$ to a solution $(A' ,\phi')$
in $\hat{{\mathcal S}}(P', E')$. By lemma \ref{l5.4} below, this solution
satisfies the condition (ii) on the intersection multiplicities.

We now use the results of section 4, namely corollary \ref{c4.3} and the
subsequent remark. These guarantee that $(A' ,\phi')$ is complex gauge
equivalent to a solution $(\tilde{A}, \tilde{\phi}) \in {\mathcal
  S}(P', E')$ of the full vortex equations. By the proof of lemma \ref{l5.7}
below, the intersection multiplicities of $\phi' (M)$ and
$\tilde{\phi}(M)$ with the submanifolds $E'_j$ are the same, and so
$(\tilde{A}, \tilde{\phi}) $ satisfies condition (ii). This proves
the existence part of the theorem.
As for the last assertion of the theorem, it follows directly from
lemma \ref{l5.7} and the fact that the $(\tilde{A}, \tilde{\phi})$ are complex
gauge equivalent to the $(A' , \phi')$. This finishes the proof for
$F_\Delta$. 

To show that the theorem remains valid for any $F$ with $\mu (F) =
\Delta$, we first remark that such an $F$ is equivariantly
biholomorphic to $F_\Delta$ \cite{Ab}. In particular, since the vortex
equations (\ref{2.0a}) and (\ref{2.0c}) only depend on the
$T^n$-action and complex structure on $F$, not on the symplectic form,
the spaces  $\hat{{\mathcal S}}(P',E')$ are the same in the $F$ and
$F_\Delta$ cases. This shows that the solution $(A' , \phi' )$
constructed above for the $F_\Delta$ case also provides a solution
for the $F$ case. Repeating the argument of the paragraph above we
conclude that the theorem also holds for $F$. 
\hfill \qed

\begin{lem}
Let $(A, \phi) \in \hat{{\mathcal S}}(P, \dot{E})$ be the pair
constructed above, and let $(A' ,\phi') \in \hat{{\mathcal S}}(P',
E')$ be its image by the map $\Upsilon$ of section 6. Then
\[
{\rm mult}_{\phi' (Z_i)} (E'_j , \phi' (M)) \ =\ a^i_j \ .
\]
\label{l5.4}
\end{lem}
\vspace{-0.4cm}
\begin{prooff}
Denote by $\pi_{\dot{E}}: \dot{E} \rightarrow E'$ the bundle map
defined in (\ref{3.3.3}), and let $E_j$ be the sub-bundle $\bigoplus_{k\ne j}
L_k$ of $E$. It follows from (\ref{6.1.4}) and the definition of $E'_j$ that
\begin{gather}
\pi_{\dot{E}}^{-1} (E'_j)\ =\ E_j \cap \dot{E} \ .
\label{6.2.2}
\end{gather}
Since the section $\phi'$ of $E'$ is by definition $\pi_{\dot{E}}
\circ \phi$, we then have that
\[
\phi'^{-1} (E'_j) \ =\ \phi^{-1} (E_j \cap \dot{E}) \ =\ \phi^{-1}
(E_j)\ .
\]
Writing this analytic hypersurface in $M$ as a union $\bigcup_{i\in I}
Z_i$ of irreducible hypersurfaces, it is tautological that
\[
E_j \cap \phi (M) \ =\ \bigcup_{i\in I}\ \phi (Z_i) \qquad {\rm and}
\qquad  E'_j \cap \phi' (M)\ =\ \bigcup_{i\in I}\ \phi' (Z_i) \ .
\]
Notice that $\phi (Z_i)$ and $\phi' (Z_i)$ are irreducible analytic
hypersurfaces in $\phi (M)$ and $\phi' (M)$, respectively, for it was
shown in the proof of lemma \ref{l4.5} that $\phi$ and $\phi'$ are
biholomorphisms onto their images.

Now, given any generic point $p \in Z_i$, let the submanifolds $E_j
\subset E$ and $E'_j \subset E'$ be locally defined around $\phi (p)$
and $\phi' (p)$ by holomorphic functions $f$ and $f'$,
respectively. This means that $f$ is a locally defined holomorphic
function whose germ at $\phi (p)$ is irreducible in the ring $\OO_{\phi (p)}
(E)$, and such that the zero locus of $f$ coincides with $E_j$ in a
neighbourhood of $\phi (p)$. Similarly for $f'$. Then from the
formulae of \cite[p. 65, 130 and 395]{G-H} it follows that
\begin{align}
{\rm mult}_{\phi (Z_i)} (E_j , \phi (M)) \ &=\ {\rm ord}_{\phi (Z_i)
  ,\, \phi (p)}\; (f|_{\phi (M)}) \ =\ {\rm ord}_{Z_i ,\, p} (\phi^\ast f)
  \qquad {\rm and}   \label{6.2.3}   \\
{\rm mult}_{\phi' (Z_i)} (E'_j , \phi' (M)) \ &=\ {\rm ord}_{\phi'
  (Z_i) ,\, 
  \phi' (p)}\; (f'|_{\phi' (M)}) \ =\ {\rm ord}_{Z_i ,\, p}
  (\phi'^\ast f')\ ,    \label{6.2.4}
\end{align}
where in the rightmost equalities we used that both $\phi$ and $\phi'$
are biholomorphisms onto their image. But it is shown in lemma \ref{l5.5}
below that if $E'_j$ is locally defined by $f'$, then $E_j$ is locally
defined by $f' \circ \pi_{\dot{E}}$. Therefore from the definition
$\phi'  = \pi_{\dot{E}} \circ \phi$ we obtain that 
\[
{\rm mult}_{\phi (Z_i)} (E_j , \phi (M)) \ = \ {\rm mult}_{\phi'
  (Z_i)} (E'_j , \phi' (M)) \ .
\] 
To finish the proof, pick local holomorphic trivializations of the
line bundles $L_j^A$ with complex coordinates $w_j$ on the
fibre. These induce a holomorphic trivialization of 
\[
E^A \ =\ L_1^A \oplus \cdots \oplus L_d^A
\]
with complex coordinates $w_1 , \ldots , w_d$ on the fibre. It is
clear that the submanifold $E_j$ is locally defined by the holomorphic
function $w_j$, and from (\ref{6.2.3}) we get that 
\[
{\rm mult}_{\phi (Z_i)} (E_j , \phi (M))\ =\ {\rm ord}_{Z_i ,\, p}\;
(\phi_j)\ =\ a^i_j \ ,
\]
where in the last equality we used that, by construction of $\phi$,
$D_j = \sum_i a^i_j \cdot Z_i$ is the divisor of the zero set of
$\phi_j$ regarded as a holomorphic section of $L_j^A$. 
\end{prooff}

\begin{lem}
In the construction above, different choices of divisors $D_j$ lead to
complex-gauge inequivalent solutions $(A', \phi') \in \hat{{\mathcal
    S}}(P',E')$. 
\label{l5.7}
\end{lem}
\begin{prooff}
Let $\{ D^{(1)}_j \}$ and  $\{ D^{(2)}_j \}$ be two sets of divisors
satisfying conditions (i) and (ii), and for $r=1,2$ let $(A_r ,
\phi_r) \in \hat{{\mathcal S}} (P_{(r)} ,E_{(r)})$ and $(A'_r , \phi'_r) \in
\hat{{\mathcal S}} (P',E')$ be the solutions obtained by the
construction above. Suppose furthermore that there exists a complex
gauge transformation $\hat{g}:M\rightarrow T^n_\CC $ such that $(A'_2
, \phi'_2) = \hat{g} (A'_1 ,\phi'_1)$.

It is shown in \cite[p. 12 and 115]{G} that there exists a subgroup $H$
of $T^d$ such that $T^d_\CC = N_\CC \times H_\CC$, and so the exact
sequence (\ref{6.1.1}) splits. In particular there exists a gauge
transformation $g:M \rightarrow T^d_\CC$ such that $\hat{g} =
\tilde{\beta}_\CC \circ g = g'$. From lemma \ref{l3.5} we then get that
$[g(A_1 ,\phi_1)]' = (A'_2 , \phi'_2)$. In particular, by lemma
\ref{l5.4} and its proof, $D^{(2)}_j$ is just the divisor of the zero set of
$g(\phi_1)_j$ regarded as a holomorphic section of
$(L_{(1)})_j^{g(A_1)}$. But it is well known that a complex gauge
transformation does not change the zero set divisor of a section of a
line bundle, and so $D^2_j$ coincides with the zero
set divisor of $(\phi_1)_j$ regarded as a holomorphic section of
$(L_{(1)})_j^{A_1} $, which by construction of $\phi_1$ is just
$D^{(1)}_j$. This proves the lemma.  
\end{prooff}

\begin{lem}
Fix a connection $A \in \AAA^{1,1} (P)$ and take the complex
structures on the bundles $E= P\times_{T^d}\CC^d$ and $E' =
P\times_{T^d} F_\Delta = P' \times_{T^d /N} F_\Delta$ to be the ones
induced by $A$, as in section 2.2.
Let $p$ be any point in $E_j \cap \dot{E}$ and suppose that, in some
neighbourhood of $\pi_{\dot{E}} (p)$, the submanifold $E'_j \subset
E'$ is locally defined by a holomorphic function $f'$. Then the
submanifold $E_j \subset E$ is locally defined by the holomorphic
function $f' \circ \pi_{\dot{E}}$ in some neighbourhood of $p$.
\label{l5.5}
\end{lem}
\begin{prooff}
Consider the projection $\pi_{\dot{E}} : \dot{E} \rightarrow E'$
defined in (\ref{3.3.3}). With respect to local trivializations of $\dot{E}$
and $E'$ induced by the same trivialization of $P$, as in (\ref{2.11}), the
map $\pi_{\dot{E}}$ is locally given by
\begin{gather}
\pi_{\dot{E}}: \UU \times \CC^d_\Delta \longrightarrow \UU \times
F_{\Delta} \ , \qquad (x,q) \mapsto (x , \pi_F (q)) \ ,
\label{6.2.5}
\end{gather}
where $\pi_F : \CC^d_\Delta \rightarrow F_\Delta$ is the holomorphic
quotient map. Now, using the local formula for the complex
structure $J(A)$ on $E$ and $E'$, using the holomorphy of $\pi_F$, and
using that $\dd \pi_F (\xi^\flat ) = \beta (\xi)^\flat$ for any vector
$\xi \in \RR^d \simeq \ttt^d$, it is not difficult to show that $\dd
\pi_{\dot{E}}$ commutes with the complex structures $J(A)$, i.e. that
$\pi_{\dot{E}}$ is a holomorphic map. From the local formula
(\ref{6.2.5}) it is 
also clear that $\pi_{\dot{E}}$ is a surjective submersion. Therefore,
by the local form of holomorphic submersions, there are neighbourhoods
$V \subset \dot{E}$ of $p$, $V' \subset E' $ of $\pi_{\dot{E}}
(p)$ and $V'' \subset \CC^{d-n}$ of the origin such that
$\pi_{\dot{E}}$ factorises as
\begin{gather}
\begin{CD}
\pi_{\dot{E}}:\ V   @>{\varphi}>> V' \times V''  
      @>>>  V' 
\end{CD} \ ,
\label{6.2.6}
\end{gather}  
where $\varphi$ is a biholomorphism and the rightmost arrow in the
canonical projection.

On the other hand, if $f'$ is a holomorphic function on $V'$ such
that
\[
V' \cap E'_j \ =\ \{ x\in V' :\  f(x)=0  \} \ ,
\]
then $f' \circ \pi_E$ is holomorphic on $V$ and, by (\ref{6.2.2}),
\[
V \cap E_j \ =\ \{ x \in V :\ f' \circ \pi_{\dot{E}} (x) =0 \}\ .
\]
Hence we only have to show that if the germ of $f'$ at $\pi_{\dot{E}}
(p)$ is irreducible, the germ of $f' \circ \pi_{\dot{E}}$ at $p$ is
irreducible as well. To do this, suppose that the germ of $f' \circ
\pi_{\dot{E}}$ is reducible. Then shrinking the neighbourhoods $V,\
V'$ and $V''$ if necessary, we may assume that 
\[
f' \circ \pi_{\dot{E}} \ =\ h_1 \cdot h_2 \qquad \quad {\rm on}\ \ V \
,
\]
where $h_1$ and $h_2$ are holomorphic functions of $V$ that vanish
at $p$. But (\ref{6.2.6}) implies that for all $q \in V'$
\[
f' (q)\ =\ f' \circ \pi_{\dot{E}} \circ \varphi^{-1} (q,0) \ =\ h_1
\circ \varphi^{-1} (q, 0)\ \cdot \ h_2 \circ \varphi^{-1} (q,0) \ .
\]
Since the two functions $h_i \circ \varphi^{-1} (\cdot ,0)$ are holomorphic
on $V'$ and vanish at $\pi_{\dot{E}}(p)$, we conclude that the germ
of $f'$ at $\pi_{\dot{E}}(p)$ is reducible as well, and this ends the proof. 
\end{prooff}

\vspace{-0.5cm}

\section{Some comments}

In this short and last section we will just make a few informal and
not completely rigorous comments about the general pattern of the
vortex solutions found in sections 3, 5.1 and 7.2.

$\ $

In all the cases the solutions are characterized by a choice of
hypersurfaces in $M$. These cannot be arbitrary hypersurfaces, but
must satisfy some topological constraints relating their Poincar\'e
duals with the Chern numbers of the bundle $P$ where the connection
$A$ is defined. Once an allowed choice of hypersurfaces is made, there
is a unique solution of the vortex equations (up to gauge equivalence)
such that the section $\phi : M \rightarrow E$ has some prescribed
values along the hypersurfaces. This prescription usually means that
along a given hypersurface the map $\phi$, which can be locally
regarded as having values on the fibre $F$, is forced to have values
on a certain complex submanifold of $F$. For $F=\CC^n$ and
$F=\CC{\mathbb P}^n$ these complex submanifolds are just the natural
$\CC^{n-1}$'s and $\CC {\mathbb P}^{n-1}$'s, respectively, contained
in $F$. When $F$ is a compact K\"ahler toric manifold these
submanifolds are the inverse images by the moment map $\mu : F
\rightarrow \RR^n$ of the $(n-1)$-dimensional faces of the Delzant
polytope $\Delta \subset \RR^n$ characterizing $F$.

The overall picture becomes clearer if one looks at simple examples,
for instance $F=\CC$ or $F= \CC {\mathbb P}^1$. In the former case
there is only one type of hypersurface to choose in $M$; along these
the Higgs field $\phi$ vanishes and they are interpreted as the
locations of the usual
vortices. In the latter case there are two types of hypersurface to 
choose in $M$: the ones taken by $\phi$ to the south-pole of $\CC
{\mathbb P}^1$ (vortices), and the ones with image the north-pole
(anti-vortices). Still in the $F= \CC {\mathbb P}^1$ case, theorem
\ref{t6.2} tells us that all hypersurfaces contribute
equally to the total energy of a solution $(A, \phi)$, independently
of their type.

$\ $

The significance of the hypersurfaces of $M$ that characterize the
solutions of the vortex equations can be better understood by varying
the real parameter $a$ in these equations. In order to do this fix the
principal bundle $P$ where the connection $A$ is defined, fix an
allowed choice of hypersurfaces, and choose a moment map $\mu : F
\rightarrow \RR^n$ such that the origin is in the interior of the
image polytope $\Delta = \mu (F)$. Then for arbitrarily large $a$ the constant
$c(a,P,M)$ of (\ref{4.1.5}) is in $\Delta$, and so solutions
exist. Furthermore, by theorem \ref{t6.2}, the energy of these
solutions tends to a finite constant as $a \rightarrow +\infty$. Now,
if the energy is to be kept constant, it is evident from (\ref{2.1})
that as $a$ grows the value of $\mu \circ \phi$ should 
approach zero almost everywhere. On the other hand we know that along
the chosen hypersurfaces $\mu \circ \phi$ has values in some face
of $\Delta$, and this is independent of $a$. Thus $\mu \circ \phi$
tends to zero everywhere except along the hypersurfaces.

Consider now the second vortex equation (\ref{2.0b}). It tells us
that, in the regions where $\mu \circ \phi \ne 0$ as $a \rightarrow
+\infty$, the quantity $\Lambda F_A$ must also become very large. Thus
in some sense the curvature of $A$, or the magnetic field, becomes
localised around the chosen hypersurfaces as $a \rightarrow
+\infty$. Notice also that, becoming localised around the
hypersurfaces, the curvature $F_A$ should be related in some way to
the Poincar\'e duals of these hypersurfaces; this is in fact what is
expressed by condition (ii) of theorems \ref{t6.1} and \ref{t5.3}.

Thus as $a$ tends to infinity the general picture is that the
solutions $(A, \phi )$ tend to the vacuum solutions of the theory ---
which are characterized by $P$ trivial, $A=0$ and $\phi = {\rm const.}
\, \in \mu^{-1}(0)$ --- except at the chosen hypersurfaces.

The opposite limit is when the parameter $a$ tends to zero. In this
case it is apparent from (\ref{2.1}) that the energy functional tends
to the pure Yang-Mills functional, and that the section $\phi$ does not
contribute to the energy. The only relevant equations are then
(\ref{2.0b}) and (\ref{2.0c}), which become the Hermite-Einstein equations.

\vskip 25pt
\noindent
{\bf Acknowledgements.}
I am pleased to thank Prof. N. S. Manton for many discussions
and some observations included in section 8. I also thank the referee
for his detailed comments. 
I am supported by \lq{\sl Funda\c{c}\~ao para a Ci\^encia e
Tecnologia}\rq, Portugal, through the research grant
SFRH/BD/4828/2001.

\newpage

\appendix

\section*{Appendices}

\section{Proof of propositions \ref{p6.4} and \ref{p6.5} }

The proof of proposition \ref{p6.4} is divided into the sequence of lemmas
\ref{l6.6}-\ref{l6.9}. Lemma \ref{l6.6} is just auxiliary; lemmas
\ref{l6.7} and \ref{l6.8} show that the map $\Upsilon$ is well defined
and a bijection, respectively; lemmas \ref{l6.10} and \ref{l6.9}
establish the second assertion of proposition \ref{p6.4}. Finally in
the last two pages of this note we prove proposition \ref{p6.5}.

The proofs of lemmas \ref{l6.6} and \ref{l6.10} will be omitted,
because they only consist of a careful unwinding of definitions. We
will not prove lemma \ref{l6.9} either, since this is a standard
result.

\begin{lem}
Let $A$  be a connection on $P$, $\phi$ a section of $E$, and denote
by $\nabla_j$ and $\nabla_{j,k}$ the connections induced by $A$ on
$V_j$ and $L_{j,k}$, respectively. Then
\[
\dd^A \phi =0 \quad \iff \quad \nabla_j (\tvp_j \circ \phi )=0 \quad
\iff \quad \nabla_{j,k} \: (\tvp_j \circ \phi )_k = 0 \quad \forall_k \ , 
\]
where $(\tvp_j \circ \phi)_k \in \Gamma (L_{j,k})$ denotes the $k$-th
component of $\tvp_j \circ \phi $ under the decomposition $V_j =
\bigoplus_{1\leq k \leq n} L_{j,k}$. Similar results hold when $(\dd^A
, \nabla_j , \nabla_{j,k})$ is substituted by $(\db^A , \nabla_j^{0,1}
, \nabla_{j,k}^{0,1})$ or $(\partial^A , \nabla_j^{1,0} ,
\nabla_{j,k}^{1,0})$. 
\label{l6.6}
\end{lem}

\begin{lem}
Let $\phi : M \rightarrow E$ be a section such that $\phi (M)
\nsubseteq E_j$ for all $j$, and let $A \in \AAA^{1,1}(P)$ be such that
$\db^A \phi = 0 $. Then $(\tilde{\varphi}_0 \circ \phi )_k \in \Gamma
(L_{0,k})$ 
is a non-zero meromorphic section of $L_{0,k}^A$ for all $k=1, \ldots , n$.
Furthermore the condition (2) in the definition of ${\mathcal C}$ is
satisfied. 
\label{l6.7}
\end{lem}
\begin{prooff}
Given a point $p\in M$ pick a neighbourhood $U$ of $p$ and
an index $j\in \{ 0, \ldots , n\}$ such that $\phi (U) \cap
E_j = \emptyset $. By lemma \ref{l6.6} we have that $(\tvp_j \circ \phi)_k$
restricted to $U$ is a holomorphic section of $L_{j,k}^A$
for all $k=1, \ldots , n$. Notice that these are non-zero sections,
otherwise $\phi (M)$ would be contained in some $E_k$, which
contradicts the assumptions. Therefore, if $j=0$, we have just proved
the result around $p$. If $j \ne 0$, then it follows from the
definition of $L_{j,k}$ (see (\ref{7.2.4})) that
$L_{0,k} = L_{j,r} \otimes L^{-1}_{j,1}$, where $r$ is appropriately
chosen among the values $k, \ k+1$ or $\emptyset$ (in which case the
term labelled by $r$ should be omitted). Furthermore, the form of the
transition functions $\varphi_0 \circ \varphi_j$ implies that
\[
(\tvp_0 \circ \phi )_k \ =\ (\tvp_0 \circ \tvp_j^{-1} \circ \tvp_j
\circ \phi)_k \ =\ (\tvp_j \circ \phi )_r \cdot [(\tvp_j \circ \phi
  )_1 ]^{-1}
\tag{A1}
\label{b2}
\]
over $U \subset \phi^{-1}(E \setminus E_j)$. But it follows from the
definitions of the induced connections that, also as
holomorphic bundles, 
\[
L_{0,k}^A \ =\ L_{j,r}^A \ \otimes\  (L^{A}_{j,1})^{-1} \, ,
\]
and hence $(\tvp_0 \circ \phi )_k$, being the quotient of two non-zero
holomorphic sections, is a non-zero meromorphic section over
$U$.

To end the proof we must now show that condition (2) is
satisfied. From formula (\ref{b2}) and the holomorphy of the $(\tvp_j
\circ \phi )_l$ it is clear that the first part of (2) is obeyed,
since the negative part of the divisor of $(\tvp_0 \circ \phi )_k$ is
just the divisor of $(\tvp_j \circ \phi )_1$. Moreover, for the
particular value $k=j$, the appropriate choice of $r$ in (\ref{b2}) is
$r=\emptyset $, and so the first factor on the r.h.s. of (\ref{b2})
should be omitted. This implies that the positive part of the divisor
of $(\tvp_0 \circ \phi )_j$ is zero, and hence the second part of
condition (2) is trivially satisfied.
\end{prooff}

\begin{lem}
Let $(\nabla_1 , \xi_1, \ldots , \nabla_n , \xi_n)$ be a multiplet
in the set ${\mathcal C}$ defined in section 5.2. Then there is a
unique connection  $A \in \AAA^{1,1} (P)$ and a unique 
section $\phi : M \rightarrow E$ such that
\begin{itemize}
\item[(i)] $\db^A \phi = 0$ and $\phi (M) \nsubseteq  E_j$ for all $j$;
\item[(ii)]$\tvp_0 \circ \phi (q) = (\xi_1 (q), \ldots , \xi_n (q)) \
  \in \ V_0 \ $ for all $\ q \in \phi^{-1} (E\setminus E_0)$;
\item[(iii)] The connection induced on $L_{0,k}$ by $A$ is $\nabla_k$.
\end{itemize}
\label{l6.8}
\end{lem}
\begin{prooff}
We start by showing how to construct the connection $A$. Let $s: U
\rightarrow P|_U $ be any local trivialization of $P$, and let
$s_k : L_{0,k} |_U \rightarrow U \times \CC$ be the associated
trivializations of $L_{0,k} = L'_k = P \times_{\rho_k} \CC$. Since the
connection $\nabla_k$ is both integrable and $h_{0,k}$-compatible, its
connection $1$-form associated to the hermitian trivialization
$s_k $ satisfies $a_k \in \Omega^1 (U ; i\RR )$ and $\dd a_k \in
\Omega^{1,1} (U ; i\RR )$. We then define the connection $A \in
\Omega^1 (P; \ttt^n )$ by the requirement
\[
s^\ast A \ =\ \left( \ldots , \sum_k (C^{-1})_{jk}\; a_k , \ldots
\right)_{1 \leq j \leq n}  \quad \in \ \Omega^1 (U ; i\RR^n \simeq
\ttt^n ).
\tag{A2}
\label{b7}
\]
The usual standard calculation can be used to show that the connection
thus defined does not depend on the trivialization $s$, and so $A$ is
globally defined on $P$. Also $F_A = \dd \, (s^\ast A)$ is in
$\Omega^{1,1}$, because $\dd a_k \in \Omega^{1,1}$, and so $A \in
\AAA^{1,1} (P)$. Finally, by formula (\ref{6.1}) the differential $\dd 
\rho_k : \ttt^n \rightarrow \ttt $ can be identified with the matrix
$(C_{kj})$ in $M_{1\times n} ({\mathbb Z})$, and so the
connection $1$-form with respect to $s_k $ of the connection
on $L_{0,k}$ induced by $A$ is 
\[
\dd \, \rho_k \circ s^\ast A \ =\ \sum_{j,l} C_{kj}
(C^{-1})_{jl}\; a_l \ =\ a_k \ .
\tag{A3}
\label{b8}
\] 
This shows that the covariant derivative induced by $A$ coincides with
$\nabla_k$, and so (iii) is satisfied. From this construction it is
also clear that the connection 1-forms (\ref{b7}) are the only ones
that satisfy (\ref{b8}) for all $k= 1, \ldots , n$. This implies that
the connection $A$ constructed above is the only connection on $P$
that satisfies (iii).

Now the construction of $\phi$. Denoting by $W_k \subset M$ the
hypersurface of singular points of the section $\xi_k$, we have that
the map
\[
\tilde{\phi} \ := \tvp_0^{-1} \circ (\xi_1 , \ldots , \xi_n)\ : \ M
\setminus \cup_k W_k \longrightarrow E\setminus E_0 
\tag{A4}
\label{b3}
\]
is smooth, and we want to show that it can be extended to a smooth
section $\phi : M \rightarrow E$.

Given an arbitrary $q \in \cup_k W_k$, choose a neighbourhood $V$ of
$q$ and holomorphic trivializations of the bundles $L_{0,k}^{\nabla_k}$
defined over $V$. With respect to these trivializations the sections
$\xi_k$ are represented by non-zero meromorphic functions $\xi'_k $,
and shrinking $V$ if necessary, these can be written as $\xi'_k =  f_k
/u_k$, where $f_k$ and $u_k$ are non-zero holomorphic functions on
$V$. Notice that, by the first part of condition (2) in the
definition of ${\mathcal C}$, the functions $u_k$ can be chosen such
that $u_1 = \cdots = u_n =: u$. Now, with respect to the
trivializations $s_k $ described at the beginning of the proof, 
the representatives of $\xi_k$ are of the form $\xi''_k = g_k \,
\xi'_k$, where $g_k$ is the transition function between $s_k$ and the
holomorphic trivialization of $L_{0,k}$. Using (\ref{7.2.9}) and the
definitions of the charts $\varphi_0$, we therefore have that 
\[
\tvp_0^{-1} \circ (\xi_1 , \ldots, \xi_n) \ =\  \big[  s, [1,
  \xi''_1 , \ldots , \xi''_n ] \; \big] \ 
 = \  \big[ s, \; [\, u, \; g_1 \, f_1 , \ldots , \; g_n \, f_n \, ]
  \; \big]
\]
over the domain $V \setminus \cup_k W_k$. But the second part of
condition (2) in the definition of ${\mathcal C}$ says that the
functions $u, f_1 , \ldots , f_n $ never vanish simultaneously, and so
the formula above shows
explicitly that $\tilde{\phi}$ can be smoothly extended to $V$. By
continuity this extension does not depend on the various choices made,
and by the arbitrariness of $q$ we actually get a global extension
$\phi :M\rightarrow E$. It is then obvious from (\ref{b3}) that property
(ii) is satisfied. Furthermore, due to their meromorphy, the
sections $\xi_k$ are zero or singular only over analytic hypersurfaces
of $M$, and so there exists $q \in M$ such that the vectors $\xi_j
(q)$ are all defined and non-zero. By formula (\ref{7.2.9}), for example, we
then get that $\phi (q) \not\in E_j$ for all $j$, and so the second
part of (i) is also satisfied. Finally, using properties (ii),
(iii) and lemma \ref{l6.6} , we conclude that $\db^A \phi = 0$ over $M
\setminus \cup_k W_k $; since both $A$ and $\phi$ are defined over the
entire $M$, by continuity we must have $\db^A \phi = 0$ on $M$, and
this establishes the existence of the section $\phi$. The uniqueness
of $\phi$ follows from condition (ii), also by continuity. 
\end{prooff}

\begin{lem}
Let $(A,\phi)$ and $(A' ,\phi')$ be two pairs in ${\mathcal B}$, and
let $(\ldots, \nabla_k , \xi_k , \ldots )$ and $(\ldots, \nabla'_k ,
\xi'_k , \ldots )$ be their images by $\Upsilon$. Then $(A,\phi)$ and
$(A',\phi')$ are complex gauge equivalent if and only if in each
$L_{0,k}$ there is a complex gauge transformation that takes
$(\nabla_k , \xi_k)$ to  $(\nabla'_k , \xi'_k)$. 
\label{l6.10}
\end{lem}

\begin{lem}
Let $L \rightarrow M$ be a complex line bundle equipped with a
hermitian metric $h$. Let also $\nabla_1$ and $\nabla_2$ be integrable
$h$-compatible connections on $L$, and $\xi_1$ and $\xi_2$ be non-zero
meromorphic sections of $L^{\nabla_1}$ and $L^{\nabla_2}$,
respectively. Then the pairs $(\nabla_1 , \xi_1)$ and  $(\nabla_2 ,
\xi_2)$ are complex gauge equivalent in $L$ if and only if the
meromorphic sections $\xi_1$ and  $\xi_2$ have the same associated
divisors in $M$.
\label{l6.9}
\end{lem}
\vspace{0.2cm}

\noindent
{\bf Proof of proposition \ref{p6.5}.} 
We consider $E$ equipped with the complex structure $J(A)$. Let
$Z\subset \phi^{-1} (E_j) \subset M$ be any irreducible analytic
hypersurface, and let $p$ be a 
generic smooth point of $Z$. Take a local complex chart $(z_1, \cdots
, z_m ,w_1, \ldots , w_n )$ of $E$, defined around $\phi (p)$, such
that the $z_k$'s are coordinates on the base $M$, the $w_k$'s are
coordinates on the fibre, and the submanifold $E_j \subset E$ is given
by the equation $w_r = 0$. We will construct such charts later
on. Since $p$ is a smooth point of $Z$, we may as well assume that $Z$
is locally defined by the equation $z_1 = 0$. With respect to this
chart, the holomorphic section $\phi :M\rightarrow E$ is locally given
by
\[
\phi : z \longmapsto (z, f_1 (z), \ldots , f_n (z)) \ ,
\]
where $z$ is the multiplet $(z_1 , \ldots , z_m)$ and the $f_k$'s are
locally defined holomorphic functions. Write
\[
f_r (z_1 , \ldots , z_m) \ =\ (z_1)^a \ g(z_1 , \ldots , z_m)\ ,
\tag{A5}   \label{b6}
\]
where $a \in {\mathbb N}_0$ and $g$ is a holomorphic function such
that $g(0, z_2 , \ldots , z_m) \not\equiv 0 $. According to the
definition of \cite[p. 130]{G-H}, we then have that
\[
a \ =\ {\rm ord}_{Z,p} (f_r)\ .
\]

The idea now is to use the formula of \cite[p. 65]{G-H} to compute the
intersection 
multiplicity  ${\rm mult}_{\phi (Z)} (E_j , \,\phi (M))$. Let $H$ be the local
submanifold of $E$ defined by the equations $z_2 = \ldots = z_m
=0$. The tangent space $T_{\phi (p)} H$ is $\CC$-generated by the
vectors
\[
\left\{ \frac{\partial}{\partial z_1}  , \frac{\partial}{\partial w_1},
\ldots , \frac{\partial}{\partial w_n}  \right\}\ .
\] 
On the other hand, since $\phi : M \rightarrow \phi (M)$ is a
biholomorphism (see the proof of lemma \ref{l4.5}), $\phi (p)$ is a
smooth point of 
the irreducible variety $\phi (Z) \subset E$, and the tangent space
$T_{\phi (p)} \phi (Z)$ is $\CC$-generated by 
\[
\left\{ \frac{\partial}{\partial z_k} + \sum_l \frac{\partial f_l}{\partial
z_k} \frac{\partial}{\partial w_l} \ :\ 2\leq k \leq m  \right\} \ .
\]
It is clear that $T_{\phi (p)}E = T_{\phi (p)}H \ +\ T_{\phi (p)}\phi
(Z)$, i.e. $H$ and $\phi (Z)$ intersect transversely at $\phi (p)$.

Now, using the inverse function theorem, it is not difficult to check
that $z_1$ together with $\tilde{w}_k = w_k - f_k (z_1 ,0,\ldots ,0)$
define a local chart for $H$ around $p$. On this chart $H\cap \phi
(M)$ is defined by the equations $\tilde{w}_1 = \cdots = \tilde{w}_n =
0$, whereas $H\cap E_j$ is defined by $w_r = \tilde{w}_r + f_r (z_1,
0, \ldots ,0) = 0$. In particular $H\cap \phi (M)$ is a submanifold of
$H$ with dimension $1$ and coordinate $z_1$. Hence, we have that for a
generic smooth point $p\in Z$,
\begin{align*}
{\rm mult}_{\phi (p)} \Big( H\cap E_j ,\ H\cap \phi (M) \Big) \ & =\ {\rm
  ord}_{\phi (p)} \Big( w_r |_{H\cap \phi (M)} \Big) \ =\ {\rm ord}_{z_1 =0}
  \Big( f_r (z_1 , 0, \ldots, 0)\Big) \\ & =\ a\ =\ {\rm ord}_{Z,p} (f_r)\ ,
\end{align*}
where we have used (\ref{b6}) and that, for a generic $p$, $g(0, \ldots, 0)
\ne 0$. It then follows from the formula of \cite[p. 65]{G-H} that
\[
{\rm mult}_{\phi (Z)} (E_j , \phi (M)) \ =\ {\rm ord}_{Z,p} (f_r)\ .
\]
We will now see what the right hand side of the above equality is.

Start by noticing that, as in the proof of lemma \ref{l6.7}, there is an
index $l \in \{0, \ldots ,n \}$ such that, around $p$, $(\tvp_l \circ
\phi )_k$ is a holomorphic section of $L^A_{l,k}$ for all
$k$. Consider then the local chart of $E$ defined by
\[
\tvp_l : E\setminus E_l \longrightarrow L^A_{l,1} \oplus \cdots \oplus
L^A_{l,n} \ ,
\]
by local holomorphic trivializations of the $L^A_{l,k}$, and by a
chart $(z_1 , \ldots ,z_m)$ of $M$ around $p$. This chart looks just
like the one described at the beginning of the proof, with $f_k$ being
the representative function of $(\tvp_l 
\circ \phi)_k$ with respect to the holomorphic trivialization of
$L^A_{l,k}$. By the form of the charts $\varphi_l$ in (\ref{7.2.8}),
the submanifold $E_j 
\subset E$ is then given by the equation $w_r =0$, where $r = j+1$ if
$j < l$ and $r=j$ if $j > l$. Notice that $j\ne l$, by the choice of
$l$. In any case, by formulae analogous to (\ref{b2}),
\[
(\tvp_l \circ \phi)_r \ =\ [(\tvp_0 \circ \phi )_l ]^{-1} \cdot
(\tvp_0 \circ \phi)_j \ =\ \xi_l^{-1} \cdot \xi_j \ ,
\]
where one should omit any $\xi$ with index $0$. Therefore
\[
{\rm ord}_{Z,p} (f_r) \ =\ {\rm ord}_Z (\tvp_l \circ \phi)_r \ =\ {\rm
  ord}_Z (\xi_j) - {\rm ord}_Z (\xi_l) \ .
\]
Finally, by the holomorphy around $p$ of 
\[
(\tvp_l \circ \phi )_k \ =\ \left\{
\begin{aligned} 
  &\xi_l^{-1} \cdot \xi_k \ \ \ \qquad {\rm if}\ \ k > l \\
  &\xi_l^{-1} \cdot \xi_{k-1} \qquad \, {\rm if}\ \; k\leq l 
\end{aligned}
\right .
\]
for all $1\leq k \leq n$, where one should omit any $\xi$ with index
$0$, we conclude that
\[
{\rm ord}_Z (\xi_l)\ =\ \min_{0\leq k\leq n} {\rm ord}_Z (\xi_k)\ .
\]
This finishes the proof.
\hfill $\qed$

\section{Complex structures from connections}

In this appendix we will prove proposition \ref{p2.2}.
This is a relatively straightforward extension of the proof given in
\cite{Ko} for the case where $F$ is a vector space.

$\ $

Use the local trivialization (\ref{2.11}) to identify $E|_U \simeq U
\times F$, where $U$ is a domain in $M$, and denote by $J_A$ the
complex structure on $U \times F$ corresponding to $J(A)$ on
$E|_U$. A vector $v \in T(U \times F)$ with components $v'$ in
$TU$ and $v''$ in $TF$ corresponds by the trivialization to the
vector
\begin{gather}
\tilde{v}\ =\ \dd\chi\; [\dd s (v') + v''] \qquad \in\ \ TE\ .
\tag{B1}   \label{a1}
\end{gather}
Now decompose the vector $\dd s (v')$ in $T(P\times F)$ as
\begin{gather}
\dd s (v')\ =\ [\dd s (v') - A(\dd s (v'))^\sharp ]\ +\ [(s^\ast
  A)(v')^\sharp - (s^\ast A)(v')^\flat ]\ +\ (s^\ast A)(v')^\flat \ ,
\tag{B2}     \label{a2}
\end{gather}
where for any $a \in \g$ we denote by $a^\sharp$ the associated
fundamental vector field on $P$, and by $a^\flat$ the vector field on
$F$ induced by the left $G$-action.
 Then the first term of (\ref{a2}) is in the horizontal space $H_A$,
the second in $\ker{\dd \chi}$, and the third in $TF$. By definition,
the complex structure $J(A)$ preserves ${\mathcal H}_A$,  $\ker{\dd
\pi_E}$, and satisfies
\[
\left\{
\begin{aligned}
& \dd \pi_E \, \circ\, J(A)\,\circ \, \dd \chi \ =\ J_M \, \circ \, \dd
\pi_E \, \circ \, \dd\chi \ =\ J_M \,\circ\,\dd\pi_P \qquad &{\rm
  on}\ \ H_A \ \ \\
& J(A)\,\circ\, \dd\chi \ =\ \dd\chi \, \circ \, J_F   \qquad &{\rm
  on}\ \ TF\ .
\end{aligned}
\right.
\]
Hence on the one hand
\begin{align*}
 \dd \pi_E \, \circ\, J(A)\,\circ \, \dd \chi\, [\; \dd s (v') - (s^\ast A)
 (v')^\sharp\; ]\ &=\ J_M (v') \\  &=\ \dd \pi_E \,\circ \,\dd \chi\,
 [\; \dd s
 \,\circ\, J_M (v') - (s^\ast A) (J_M v')^\sharp\; ] , 
\end{align*}
and since $\dd \pi_E : {\mathcal H}_A \rightarrow TM$ is an
isomorphism,
\begin{align}
J(A)\,\circ \, \dd \chi\; [\; \dd s (v') - (s^\ast A)(v')^\sharp\; ]\ &=\
 \dd \chi\; [\, \dd s \,\circ\, J_M (v') - (s^\ast A) (J_M v')^\sharp\, ]
 \nonumber  \\
&=\ \dd \chi\; [\,\dd s\,\circ\, J_M (v') - (s^\ast A) (J_M v')^\flat\, ]\ .
\tag{B3}    \label{a3}
\end{align}
On the other hand
\begin{gather}
J(A)\,\circ \, \dd \chi\; [\, v'' + (s^\ast A) (v')^\flat \, ]\ =\ \dd\chi
\,\circ\,  J_F \; [\, v'' + (s^\ast A) (v')^\flat\, ] \ .
\tag {B4}    \label{a4}
\end{gather}
Using (\ref{a1})-(\ref{a4}) and the convention (\ref{3.1.0}) we then
obtain that 
\begin{align*}
[J_A(v)]^\sim      \ &=\ J(A)\,[\tilde{v}]\\ &=\ \dd\chi \,
\left\{ \dd s 
\, \circ\, J_M (v') + J_F \, [v'' + (\,(s^\ast A) (v') + i (s^\ast
  A)(J_M v')\, )^\flat\, ] \right\} \\
&=\ \left\{ J_M (v') + J_F\,[v'' + 2 (s^\ast A)^{0,1} (v')^\flat\, ]
\right\}^\sim   \ .
\end{align*}
Since the map $v\mapsto \tilde{v}$ of (\ref{a1}) is an isomorphism, we
finally get that
\begin{gather}
J_A \; (v) \ =\ J_{M\times F} \, [v + 2 (s^\ast A)^{0,1} (v)^\flat ] \qquad
\in \ \ T(U \times F)\ ,
\tag{B5}      \label{a5}
\end{gather}
where $J_{M\times F}$ denotes the natural complex structure on the
product $M\times F$. This is an explicit formula for the complex
structure $J(A)$ as seen through the trivialization $E|_U \simeq U
\times F$; we will rely on it to study the integrability of $J(A)$.

$\ $

Let $U' \times V' \subset U \times F$ be the domain of a complex
chart $\{ z^j , w^r \}$ of $U \times F$, where the $z^j$'s and the
$w^r$'s are coordinates on $U'$ and $V'$, respectively. For a
given basis $\{ \xi_k \}$ of the Lie algebra $\g$, write over $V'$
\[
(\xi_k)^\flat \ =\ \alpha_k^r \, \frac{\partial}{\partial w^r}\ .
\] 
Since the action transformations $\rho_g : F \rightarrow F$ are
holomorphic, the vector fields $\xi_k^\flat $ and the functions
$\alpha_k^r$ are both holomorphic over $V'$. Now consider the set of
1-forms over $U' \times V'$
\[
D\ =\ \left\{ \dd\,z^j,\ \dd\, w^r + (s^\ast A^k)^{0,1}\, \alpha_k^r \
:\ j=1,\ldots , \dim_\CC{M}\ \ {\rm and}\ \ r=1, \ldots ,\dim_\CC{F}
\right\}\ .
\]
Using formula (\ref{a5}) and the fact that, with respect to $J_{M\times F}$,
the forms $\dd z^j$ and $\dd w^r$ are of type $(1,0)$ while $(s^\ast
A)^{0,1}$ is of type $(0,1)$, it is not difficult to check that all
the forms in $D$ are of type $(1,0)$ with respect to $J_A$. It is also
apparent that they are linearly independent at each point of $U'
\times V'$, and so we conclude that they generate the submodule
$\Omega_A^{1,0} ( U' \times V' )$ of $\Omega^1 ( U' \times V')$,
where the splitting $\Omega^1 = \Omega_A^{1,0} \oplus
\Omega_A^{0,1}$ is with respect to the complex structure $J_A$. Now by
basic results on complex manifolds, the complex structure $J_A$ is
integrable iff $\dd\,\Omega_A^{1,0}\subseteq \Omega_A^{2,0} \oplus
\Omega_A^{1,1}$, and this last condition is clearly
equivalent to 
\begin{gather}
\dd\,\eta \ \in\ \Omega_A^{2,0} \oplus \Omega_A^{1,1} \qquad \ {\rm
  for\ all}\ \eta \in D \ .
\tag{B6}       \label{a6}
\end{gather}
Thus the problem is reduced to a simple computation. This computation
yields
\[
\dd\, (\dd\, z^j)\ =\ 0
\]
and, modulo $\Omega_A^{2,0} \oplus \Omega_A^{1,1}$,
\begin{align*}
 \dd\, (\dd w^r + (s^\ast A^k)^{0,1} \alpha_k^r \, ) \ &=
\ \frac{\partial \alpha_k^r}{\partial w^l}\, (\dd w^l + (s^\ast
A^n)^{0,1} \alpha_n^l\, )\, \wedge \, (s^\ast A^k)^{0,1} \ - \\
&  \  -\, 
\frac{1}{2}\, \left( \frac{\partial \alpha_k^r}{\partial w^l} \,
\alpha_n^l - \frac{\partial \alpha_n^r}{\partial w^l}\, \alpha_k^l
\right) \, (s^\ast A^n)^{0,1} \wedge  (s^\ast A^k)^{0,1}\, + \,
\alpha^r_k \,\dd (s^\ast A^k)^{0,1}   \\
&=\ -\, \frac{1}{2}\, [\xi_n^\flat , \xi_k^\flat ]^r \, (s^\ast A^n)^{0,1}
\wedge  (s^\ast A^k)^{0,1} \, + \, (\xi_k^\flat )^r \, \db (s^\ast
A^k)^{0,1}   \\
&=\ \left\{ \frac{1}{2}\, (\,[\xi_n , \xi_k ]^\flat\, )^r \, (s^\ast
A^n)\, \wedge \, (s^\ast A^k) \, +\, (\xi_k^\flat\, )^r \, \dd (s^\ast
A^k) \right\}^{0,2}    \\
&=\ \left( F_A^{0,2}  \right)^k \, \left( \xi_k^\flat  \right)^r  \ ,
\end{align*}
where $F_A$ is the curvature form on the base $M$. Hence (\ref{a6}) is
satisfied if and only if 
\begin{gather}
(F_A^{0,2})^k \xi_k^\flat \  =\  0 \qquad {\rm on}\ \ U' \ ,
\tag{B7}        \label{a7}
\end{gather}
and the first statement of the proposition is
obviously true. As for the second statement, suppose that at least one
point in $F$ has a discrete isotropy group. Then by well known results
almost every point in $F$ has this property \cite[p. 179]{Bred}. In particular,
for any point $p$ in an open dense subset of $M$, the linear map $\g
\rightarrow T_p F$ defined by $\xi \mapsto \xi^\flat$ is injective and
the vectors $\xi_k^\flat |_p$ are linearly independent. In this case
if $J(A)$ is integrable, i.e. if (\ref{a7}) is satisfied, necessarily
$F_A^{0,2} =0$.

\end{document}